\documentclass[11pt]{article}
\usepackage{amsfonts}

\usepackage{graphics}

\usepackage{indentfirst}
\usepackage{cite}
\usepackage{latexsym}
\usepackage{amsmath}
\usepackage{amssymb}
\usepackage{showkeys} 
\usepackage[dvips]{epsfig}
\usepackage{amscd}
\def\on{\bar\rho}
\newtheorem{theorem}{Theorem}[section]
\newtheorem{remark}{Remark}[section]

\newtheorem{definition}{Definition}[section]
\newtheorem{lemma}[theorem]{Lemma}

\newtheorem{proposition}[theorem]{Proposition}

\newcommand{\n}{\rho}

\newcommand{\ti}{\tilde}
\newcommand{\mr}{\mathbb{R}}
\newcommand{\lm}{\lambda}

\renewcommand{\div}{ {\rm div }  }

\newcommand{\pa}{\partial}
\renewcommand{\r}{\mathbb{R}}

\renewcommand{\b}{B_{N_*}}

\newcommand{\ia}{\int_0^T}

\newcommand{\bt}{\begin{theorem}}
\newcommand{\bl}{\begin{lemma}}
\newcommand{\el}{\end{lemma}}
\newcommand{\et}{\end{theorem}}
\newcommand{\ga}{\gamma}

\newcommand{\al}{\alpha}
\newcommand{\de}{\delta}
\newcommand{\ve}{\varepsilon}
\newcommand{\la}{\label}

\newcommand{\ol}{\overline}

\newcommand{\bn}{\begin{eqnarray}}
\newcommand{\en}{\end{eqnarray}}
\newcommand{\bnn}{\begin{eqnarray*}}
\newcommand{\enn}{\end{eqnarray*}}

\newcommand{\bnnn}{\begin{eqnarray*}}
\newcommand{\ennn}{\end{eqnarray*}}

\newcommand{\ba}{\begin{aligned}}
\newcommand{\ea}{\end{aligned}}
\newcommand{\be}{\begin{equation}}
\newcommand{\ee}{\end{equation}}
\def\O{{\r^2 }}
\def\p{\partial}
\def\norm[#1]#2{\|#2\|_{#1}}

\def\o{\omega}

\newcommand{\no}{\nonumber\\}

\newcommand{\si}{\sigma}

\def\la{\label}

\def\na{\nabla}
\def\on{\bar\n}

\makeatletter      
\@addtoreset{equation}{section}
\makeatother

\makeatletter      
\@addtoreset{equation}{section}
\makeatother       

\title{Global existence and large time asymptotic behavior of strong solutions to the 2-D compressible  magnetohydrodynamic equations with  vacuum \thanks{This work was partially supported by
  NNSFC (Grant Nos.  11171228 \& 11301431),
  the Fundamental Research Funds for the Central Universities (Grant Nos. 2012121005).}}

\author{Boqiang L{\small V}\thanks{College of Mathematics and Information
Science, Nanchang Hangkong University, Nanchang, 330063, China({\tt lvbq86@163.com}).
 }
\quad Xiaoding S{\small HI}\thanks{Department of Mathematics, School of Science, Beijing University of Chemical Technology, Beijing 100029, China ({\tt shixd@mail.buct.edu.cn})}
\quad Xinying X{\small U}\thanks{Corresponding author, School of Mathematical Sciences, Xiamen University, Xiamen, China ({\tt xinyingxu@xmu.edu.cn}) }
 }

\date{ }

\begin{document}
\maketitle

 \begin{abstract} This paper concerns  the Cauchy problem of the magnetohydrodynamic equations for viscous compressible barotropic flows in two or three spatial dimensions with vacuum as far field density. For two spatial dimensions, we establish the global existence
and uniqueness of strong solutions  (which may be of possibly large oscillations)
  provided the smooth initial data are of small total energy, and obtain some a
priori decay with rates (in large time) for  the pressure,  the spatial gradient
of both the velocity field and the magnetic field. Moreover, for three spatial dimensions case, some  decay rates  are also obtained.
 \end{abstract}


\section{Introduction}

We consider the  magnetohydrodynamic (MHD) equations
\be \la{a1}  \begin{cases} \n_t + \div(\n u) = 0,\\
 (\n u)_t + \div(\n u\otimes u) + \nabla P(\n) = \mu\Delta u + (\mu + \lambda)\nabla(\div u)+(\na\times H)\times H,\\
 H_t-\na \times(u\times H)=-\na \times(\nu \na \times H),\quad \div H=0,
\end{cases}\ee for viscous  compressible magnetohydrodynamics flows. Here, $t\ge 0$ is time,   $x\in\r^2$  is the spatial coordinate,  and $\n=\n(x,t),$
$u=(u^1, u^2)(x,t),  $ $H=(H^1,H^2)(x,t)$,  and \be P(\n)=R\n^\ga\,\,( R>0, \ga>1)  \ee   are the fluid
density, velocity, magnetic field and pressure, respectively.  Without loss of generality, we assumed that $R=1.$ The constant
viscosity coefficients $\mu$ and $\lambda$ satisfy the physical
restrictions: \be\la{h3} \mu>0,\quad \mu +\lambda\ge 0.
\ee
The constant $\nu>0$ is the resistivity coefficient which is inversely proportional to the electrical conductivity constant and acts as the magnetic diffusivity of magnetic fields.
We consider  the Cauchy problem
for (\ref{a1}) with $(\n,u,H)$  vanishing at infinity (in some weak sense) with given initial data $\n_0$, $u_0,$  and $H_0$ as
\be
\la{h2} \n(x,0)=\n_0(x), \quad \n u(x,0)=\n_0u_0(x),\quad H(x,0)=H_0, \quad x\in  \r^2.
\ee

There have been huge literatures on the compressible MHD problem \eqref{a1} by many physicists and mathematicians
due to its physical importance, complexity, rich phenomena and mathematical challenges, see for
example, \cite{ca,df2006,cw,cw2,fjn,fy1,fy2,hw1,hw2,ko,ka,lyz,chl,uks,vk,w,zjx,zz} and the references therein. In particular, if there is no electromagnetic effect,
i.e., $H=0$, then \eqref{a1} reduces to the compressible Navier-Stokes equations for barotropic flows,
which have also been discussed in numerous studies, see for example, \cite{choe2003,F2,F1,H3,hof2002,hx2,hlx,hlx1,hlw,hl,lx1,lx,liliang,lxz,L2,L1,M1,xin98,xin13} and the references therein.
The issues of well-posedness and dynamical behaviors of MHD system are rather
complicated to investigate because of the strong coupling and interplay interaction
between the fluid motion and the magnetic field.  Now, we briefly recall some results
concerned with the multi-dimensional compressible MHD equations which are more relatively with our problem.
The local strong solutions to the compressible MHD with
large initial data were obtained, by Vol'pert-Khudiaev \cite{vk} as the initial density
is strictly positive and by Fan-Yu \cite{fy2} as the initial density may contain vacuum,
respectively. And   recently, the local existence  of strong and classical solutions to the two-dimensional  compressible MHD equations with vacuum as far field density has been studied in \cite{chl}.  The global existence
of  solutions to the compressible MHD equations were obtained in many works:
 Kawashima \cite{ka} firstly obtained the global existence when the initial data are close to a non-vacuum equilibrium in $H^3$-norm; Hu-Wang \cite{hw1,hw2} and Fan-Yu \cite{fy1} proved the global existence
of renormalized solutions  under the  general large initial
data assumptions; For the case that the initial density is allowed to vanish and even has compact support,
Li-Xu-Zhang \cite{lxz} established the  global existence and uniqueness of
classical solutions with constant state as far field which could be either vacuum or nonvacuum
to (1.1)-(1.4) in three-dimensional space with smooth initial data which are of
small total energy but possibly large oscillations, which generalized the results of Huang-Li-Xin \cite{hlx1} for barotropic  compressible Navier-Stokes equations to the compressible MHD ones.
Moreover, it was also showed in \cite{lxz} that for any $p > 2$, the following   large-time behavior of the solution holds:
\be\la{intro}
 \lim_{t\rightarrow\infty}\left(\|P(\n)-P( \tilde{\n})\|_{L^p(\mathbb{R}^3)}+\|\na u\|_{L^2(\mathbb{R}^3)}+\|\na  H\|_{L^2(\mathbb{R}^3)}\right)=0
\ee
where $\tilde{\n}$ is the constant far field density.

  For two-dimensional
problems, only in the case that the far field  density is away from vacuum, the techniques of \cite{lxz} can be modified directly since at this case, for any $p\in [2,\infty),$ the $L^p$-norm of a function $u$ can be bounded   by $\|\n^{1/2} u\|_{L^2}$  and $\|\na u\|_{L^2},$    and the similar results can be obtained. However, when the far field density is vacuum, it seems difficult to  bound the $L^p$-norm of $u$ by $\|\n^{1/2} u\|_{L^2}$  and $\|\na u\|_{L^2} $ for any $p\ge 1,$ so the global existence and large time behavior of strong or classical solutions to the Cauchy problem are much more subtle and remain open. Therefore,  the main  aim  of this paper is to  study the  global  existence and large time behavior of   strong solutions to   \eqref{a1}-\eqref{h2}  in some homogeneous  Sobolev spaces in two-dimensional space with vacuum as far field density. Although recently,  for the  two-dimensional Cauchy problem of barotropic compressible Navier-Stokes equations with vacuum as far field density, Li-Xin \cite{lx1} obtained both the global existence of strong solutions  and  the decay rates of the pressure and the gradient of velocity. However, their theory cannot be applied directly to the MHD ones.

Before stating the main results, we first explain the notations and
conventions used throughout this paper. For $R>0$, set
$$B_R  \triangleq\left.\left\{x\in\r^2\right|
\,|x|<R \right\} , \quad \int fdx\triangleq\int_{\r^2}fdx.$$ Moreover, for $1\le r\le \infty, k\ge 1, $ and $\beta>0,$ the standard homogeneous and inhomogeneous Sobolev spaces are defined as follows:
   \bnn  \begin{cases}L^r=L^r(\r^2 ),\quad D^{k,r}=D^{k,r}(\r^2)=\{v\in L^1_{\rm loc}(\r^2)| \na^k v\in L^r(\r^2)\}, \\ D^1=D^{1,2},\quad
W^{k,r}  = W^{k,r}(\r^2) , \quad H^k = W^{k,2} , \\ \dot H^\beta=\left\{f:\r^2
 \rightarrow \r\left|\|f\|^2_{\dot H^\beta}=
 \displaystyle{\int} |\xi|^{2\beta}|\hat f(\xi)|^2d\xi<\infty\right.
 \right\} ,\end{cases}\enn  where $\hat f$ is the Fourier transform
 of $f.$  Next, we give the definition of strong solution to \eqref{a1} as follows:
\begin{definition} If  all derivatives involved in \eqref{a1} for $(\rho,u,H)  $  are regular distributions, and   equations  \eqref{a1} hold   almost everywhere   in $\r^2\times (0,T),$ then $(\n,u,H)$  is called a  strong solution to  \eqref{a1}.
\end{definition}

The initial total energy  is defined as: \bnn  C_0 =
\int_{\r^2}\left(\frac{1}{2}\n_0|u_0|^2 +\frac{1}{2}|H_0|^2+ \frac{1}{\ga-1}P(\n_0) \right)dx. \enn Without loss of generality, assume that the initial density $\n_0$ satisfies
\be\la{oy3.7} \int_{\r^2} \n_0dx=1,\ee  which implies that there exists a positive constant $N_0$ such that  \be\la{oy3.8} \int_{B_{N_0}}  \n_0  dx\ge \frac12\int\n_0dx=\frac12.\ee

We can now state our main result in this paper, concerning the global existence  of  strong solutions.
\begin{theorem}\la{th1} In addition to \eqref{oy3.7} and  \eqref{oy3.8},  suppose that  the initial data $(\n_0,u_0,H_0)$ satisfy  for any given numbers $M>0,$ $\on\ge 1,$   $a>1,~q>2$,  and $\beta\in (1/2,1],$
\begin{equation} \la{co1}\left\{\begin{array}{lll} 0\le   \n_0\le\bar{\n}, \quad  \bar x^a \rho_0\in   L^1 \cap H^1\cap W^{1,q},\\(u_0,H_0)\in \dot H^\beta \cap D^1,\quad \n_0^{1/2}u_0  \in L^2,\quad \bar x^{{a}/{2}}H_0 \in L^2,
\end{array}\right.
\end{equation}
 and that \be\la{h7}
   \|u_0\|_{\dot H^\beta}+\|H_0\|_{\dot H^\beta}+ \|\n_0\bar x^a\|_{L^1} + \||H_0|^2\bar x^a\|_{L^1}\le M ,  \ee   where
  \be\la{2.07} \bar x\triangleq(e+|x|^2)^{1/2} \log^2 (e+|x|^2) .\ee Then there exists a positive constant $\ve$ depending
 on  $\mu ,  \lambda ,   \ga ,  a ,  \nu,  \on, \beta, N_0,$ and $M$  such that if
 \be\la{i7}
     C_0\le\ve,
   \ee  the problem  \eqref{a1}--\eqref{h2} has a unique global strong solution $(\n,u,H)$ satisfying   for any $0<   T<\infty,$
\be\la{h8}
0\le\n(x,t)\le 2\bar{\n},\quad  (x,t)\in \O\times[0,T],
\ee
\be\la{1.10}
\begin{cases}
\rho\in C([0,T];L^1 \cap H^1\cap W^{1,q} ),\\
\bar x^a\rho\in L^\infty( 0,T ;L^1\cap H^1\cap W^{1,q} ),\\
\sqrt{\n } u,\,\na u,\, \bar x^{-1}u,\,    \sqrt{t} \sqrt{\n}  u_t \in L^\infty(0,T; L^2 ), \\
H, H^2, H \bar{x}^{a/2}, \na H,  \sqrt{t}H_t\in L^\infty( 0,T ;L^2), \\
\na u\in  L^2(0,T;H^1)\cap  L^{(q+1)/q}(0,T; W^{1,q}), \\
\na H\in L^2(0, T; H^1),\\
\sqrt{t}\na u\in L^2(0,T; W^{1,q} ),  \\
\sqrt{\n} u_t, \,\na H\bar{x}^{a/2}, \,  \sqrt{t}\na u_t ,\, \sqrt{t}\na H_t ,\,  \sqrt{t} \bar x^{-1}u_t\in L^2(\O\times(0,T)),\\
   \end{cases}\ee
and
\be\la{l1.2}
\inf\limits_{0\le t\le T}\int_{B_{N_1(1+t)\log^\al(e+t) }}\n(x,t) dx\ge \frac14,
\ee
for any $\al>1$  and  some positive constant $N_1$ depending only on $\al,N_0, $ and $M.$ Moreover, $(\n,u, H)$ has the following decay
rates,  that is,  for $t\ge 1,$
 \be \la{lv1.2}
\begin{cases}
\|\na H(\cdot,t)\|_{L^2}\le Ct^{-1/2},\\
\|\na u(\cdot,t)\|_{L^p}  \le C(p)t^{-1+1/p},\quad  for~p\in[2,\infty),\\
\|P(\cdot,t)\|_{L^r}\le C(r)t^{-1+1/r},\quad  for~r\in(1,\infty),\\
\|\na \omega(\cdot,t)\|_{L^2}+\|\na F(\cdot,t)\|_{L^2}\le Ct^{-1},
\end{cases} \ee
where \be \la{hj1} \o\triangleq  \pa_1u^2-\pa_2u^1 ,\quad F\triangleq(2\mu+\lambda)\div u-P-\frac{1}{2}|H|^2,\ee are respectively  the vorticity and the so-called effective viscous flux, and $C$ depends on  $  \mu ,  \lambda ,\nu,   \ga ,  a ,  \on, \beta, $ $ N_0,$ and $M.$
\end{theorem}

For the three-dimensional case, that is,  $\Omega=\r^3,$   we have the following results concerning the decay properties of the global classical solutions whose existence can be found in \cite{lxz}.
\begin{theorem}\la{thv}  Let $\Omega=\r^3.$
  For given numbers $M>0,$ $\on\ge 1,$
 $\beta\in (1/2,1],$ and $q\in (3,6),$ suppose that
    the initial data $(\n_0,u_0,H_0)$ satisfy
\be \la{cvo1}  \rho_0 ,\, P(\n_0) \in H^2\cap W^{2,q},\quad P(\n_0),\,\n_0 | u_0|^2  \in L^1,\quad   u_0, H_0 \in \dot
H^\beta,\quad \nabla u_0,\nabla H_0 \in H^1 , \ee
 \be\la{hv7} 0\le \rho_0\le \bar{\rho},\quad
   \|u_0\|_{\dot H^\beta}+\|H_0\|_{\dot H^\beta} \le M,   \ee
and the compatibility condition \be\la{cvo2} -\mu\triangle u_0
-(\mu+\lambda)\nabla \div u_0-(\na\times H_0)\times H_0 + \nabla P(\rho_0) = \rho^{1/2}_0g, \ee
for some $ g\in L^2.$ Moreover, if $\ga>3/2,$ assume that \be \n_0\in L^1.\ee
 Then there exists a positive constant $\ve$ depending
 on $\mu,\lambda,   \ga, \nu, $  $\on, \beta,$ and $M$  such that if
 \be\la{iv7}
     C_0\le\ve,
   \ee  the Cauchy problem
  (\ref{a1})-(\ref{h2})
  has a unique global  classical  solution $(\rho,u,H)$ in
   $\r^3\times(0,\infty)$ satisfying for
  any $0<\tau<T<\infty,$
  \be\la{hv8}
  0\le\rho(x,t)\le 2\bar{\rho},\quad x\in \r^3,\, t\ge 0,
  \ee
   \be
   \la{hv9}\begin{cases}
   \rho  \in C([0,T];L^{3/2}\cap H^2\cap W^{2,q}),\\ P\in C([0,T];  L^1\cap H^2\cap W^{2,q}) , \quad  u\in L^\infty(0,T;L^6),\\
\na  u\in L^\infty(0,T;H^1)\cap L^2(0,T;H^2)\cap    L^\infty(\tau,T;H^2\cap W^{2,q}),\\
\na  u_t\in L^2(0,T;L^2) \cap L^{\infty}(\tau,T;H^1)\cap H^1(\tau,T;L^2),\\
H\in C([0, T]; H^2)\cap L^\infty(\tau, T; H^3),\\
H_t\in C([0, T]; L^2)\cap H^1(\tau, T; L^2)  .\end{cases} \ee
 Moreover, for $r\in (1,\infty),$ there exist positive constants $C(r)$ and $C$ depending   on $\mu,\lambda,   \ga, $  $\on, \beta,$ and $M$ such that for $t\ge 1,$
\be \la{lvy8}\begin{cases}
\|\na H(\cdot,t)\|_{L^p}  \le C t^{-1+(6-p)/(4p)}, \mbox{ for }   p\in [2,6] , \\
\|\na u(\cdot,t)\|_{L^p}  \le C t^{-1+1/p}, \mbox{ for }   p\in [2,6] , \\
\|P(\cdot,t)\|_{L^r}\le C(r)t^{-1+1/r},  \mbox{ for }   r\in (1,\infty),\\
\|\nabla F(\cdot,t)\|_{L^2}+ \|\nabla \omega(\cdot,t)\|_{L^2}\le C t^{-1},
\end{cases} \ee
 where $F,\omega$ defined in \eqref{hj1}, if $\ga>3/2,$ $C(r)$ and $C$ both depend on $\|\n_0\|_{L^1(\r^3)}$ also.
\end{theorem}

\begin{remark} \la{re1}
When $H=0$, i.e., there is no electromagnetic field effect, \eqref{a1} turns to be the  compressible Navier-Stokes equations, and Theorems \ref{th1} and \ref{thv} are the same as those results of   Li-Xin \cite{lx1}. Roughly speaking, we generalize  the results of\cite{lx1} to the compressible MHD equations.
\end{remark}

\begin{remark}\la{re2}
It should be noted here that the large time decay rate estimates \eqref{lv1.2}   plays a crucial role in deriving the global existence of strong
solutions to the two-dimensional problem (1.1)-(1.4), which is  completely different from  the three-dimensional case (\cite{lxz}). More precisely,    the global existence of classical solutions to (1.1)-(1.4) in \cite{lxz}
was achieved without any bounds on the decay rates of the solutions partially due to the
a priori $L^6$-bounds on the velocity field and the magnetic field.\end{remark}

\begin{remark}\la{re3}
We should point  out that the large time asymptotic decay with rates
of the global strong   solutions, \eqref{lv1.2} and \eqref{lvy8},  yield in particular that the $L^2$-norm of  the pressure, the gradient
of the velocity and the magnetic  decay in time with a rate $t^{-1/2}$, and the gradient of the vorticity and
the effective viscous flux decay faster than themselves. As will be seen in the proof,
the  large time asymptotic decay are mainly controlled by the  decay rate of the $L^2$-norm of  the pressure.
When $H=0$, the large time  decay rates \eqref{lv1.2} and \eqref{lvy8} are the same as theirs in \cite{lx1}. However, the decay rates of the magnetic field for large time in  \eqref{lv1.2} and \eqref{lvy8}  are completely new for both  the two and three spatial dimensions compressible MHD equations.
\end{remark}

\begin{remark}\la{re4}
Similar as \cite{chl}, if the initial data $(\n_0,u_0,H_0)$ satisfy some additional regularity and compatibility conditions,  the global strong solutions obtained by Theorem \ref{th1} become classical ones.
\end{remark}

We now make some comments on the analysis of this paper. Note that for initial data
in the class satisfying \eqref{co1} and \eqref{h7} except $(u_0, H_0)\in \dot H^\beta $, the local existence
and uniqueness of strong and classical   solutions to the Cauchy problem, (1.1)-(1.4), have been
established recently in \cite{chl}. To extend the strong solution globally in time, one
needs some global a priori estimates on smooth solutions to (1.1)-(1.4) in suitable higher
norms. It turns out that  the key issue here is to
derive both the time-independent upper bound for the density and the time-depending
higher norm estimates of the smooth solution $(\n, u, H).$ To this end, on the one hand, we try to   adapt   some basic ideas used in \cite{lxz}. However, new difficulties arise  in the two-dimensional case, since the analysis in \cite{lxz} relies heavily on the basic fact that,  the
$L^6$-norm of $v\in D^1(\mathbb{R}^3)$ can be bounded by the $L^2$-norm of the gradient of $v$ which fails
for $v\in D^1(\mathbb{R}^2)$. On the other hand, compared with  the two-dimensional   compressible Navier-Stokes equations considered by Li-Xin (\cite{lx1}),  for the compressible MHD equations, the strong
coupling between the velocity vector field and the magnetic field such as $\na\times(u\times H)$    and  $(\na \times  H)\times H,$  will bring out some new difficulties. Therefore, motivated by
  \cite{lx1},  we try to obtain that the $L^2$-norm
in both space and time of the pressure is time-independent (see  \eqref{h27}). However,
 the usually $L^2$-norm (in both space and time) of $H_t$    cannot be directly estimated due to the  strong
coupled term  between the velocity vector field and the magnetic field,   $\na\times(u\times H).$
The key observation to overcome this difficulty is as follows: Instead of  estimating  the $L^2$-norm of $H_t,$   we multiply  the magnetic equations  by $\Delta H$ and $H\Delta |H|^2$ respectively (see \eqref{mm1} and \eqref{nlv4}),  and succeed in  controlling   the coupled term $\na\times(u\times H)$  by the gradient of both the velocity and the magnetic after  integration by parts.  This yields  some new desired a priori estimates of the $L^2$-norm of $|H||\Delta H|$ in both space and time. In fact, this is the first key observation of this paper.
 Next, our second key point  is to get the $H^1$-norm of the effective viscous flux decays at the rate of $t^{-1/2}$ for large time (see \eqref{hg2}).
This is completed by  driving the rates of decay for not only  $\na u$ and $P$ (compared with \cite[Lemma 3.4]{lx1}) but also $H$ and $\na H$.
Indeed, we prove that the $L^2$-norm of $|H|^2$ and $|H||\na H|$ decay  at the rates of  $t^{-1/2}$ and  $t^{-1}$ respectively  (see \eqref{ly8}).
 Then, using the expansion rates of the essential support of the density (see \eqref{uq2} or \cite[(3.39)]{lx1})
for large time, we obtain the bound of  the $L^p$-norm of the gradient of the effective viscous flux (see \eqref{z.2}).
 Based on these key ingredients,
we are able to obtain the  estimates on $L^1(0,\min\{1, T\}; L^\infty(\mathbb{R}^2))$-norm and the
time-independent ones on $L^4( \min\{1, T\},T; L^\infty(\mathbb{R}^2))$-norm of the effective viscous flux
(see \eqref{hg1}). Then, motivated by \cite{lx}, with the help of  these estimates and Zlotnik's
inequality (see Lemma \ref{le1}),  we obtain  the desired time-uniform upper bound of  the density, which is the key for global estimates of strong solutions.
 Then,   similar to arguments as \cite{hx2,hlx,hlx1,lx1,lxz}, the next main step  is to bound
the gradients of the density,  velocity and magnetic.
More precisely, such bounds can be obtained by solving a logarithm Gronwall's inequality based on a Beale-Kato-Majda
type inequality (see Lemma \ref{le9}) and the a priori estimates we have just derived, and
moreover, such a derivation yields simultaneously also the bound for
 $L^1(0,T; L^\infty(\mathbb{R}^2))$-norm of the gradient of the velocity, see Lemma \ref{le4} and its proof.
 Finally, our third new observation of this paper is to obtain the $L^2$-norm of $\bar{x}^{a/2}H$ and $\bar{x}^{a/2}\na H$ (see \eqref{gj10} and \eqref{gj10'}), which will be used in the estimates of the  $H^1$-norm of the  gradient of both the velocity and the magnetic (see \eqref{gj13'}).

The rest of the paper is organized as follows: In Section \ref{se2}, we
collect some elementary facts and inequalities which will be needed
in later   analysis. Sections \ref{se3} and \ref{se5} are devoted to deriving the necessary
a priori estimates on strong solutions which are needed to extend
the local solution to all time. Then finally, the main results,
Theorems \ref{th1}-\ref{thv}, are proved in Section \ref{se4}.

\section{Preliminaries}\la{se2}

In this section, we will recall some  known facts and elementary
inequalities which will be used frequently later.

We begin with the local existence of strong and classical solutions whose proof   can be found in \cite{chl}.

\begin{lemma}   \la{th0}  Assume  that
 $(\n_0,u_0,H_0)$ satisfies \eqref{co1}  except $(u_0, H_0)\in \dot H^\beta.$   Then there exist  a small time
$T >0$    and a unique strong solution $(\rho , u,H )$ to the
  problem   \eqref{a1}-\eqref{h2}  in
$\O\times(0,T )$ satisfying \eqref{1.10} and  \eqref{l1.2}. 
 \end{lemma}

Next, the following well-known Gagliardo-Nirenberg inequality (see \cite{nir})
  will be used later.

\begin{lemma}
[Gagliardo-Nirenberg]\la{l1} For  $p\in [2,\infty),q\in(1,\infty), $ and
$ r\in  (2,\infty),$ there exists some generic
 constant
$C>0$ which may depend  on $p,q, $ and $r$ such that for   $f\in H^1({\O }) $
and $g\in  L^q(\O )\cap D^{1,r}(\O), $    we have \be
\la{g1}\|f\|_{L^p(\O)}^p\le C \|f\|_{L^2(\O)}^{2}\|\na
f\|_{L^2(\O)}^{p-2} ,\ee  \be
\la{g2}\|g\|_{C\left(\ol{\O }\right)} \le C
\|g\|_{L^q(\O)}^{q(r-2)/(2r+q(r-2))}\|\na g\|_{L^r(\O)}^{2r/(2r+q(r-2))} .
\ee
\end{lemma}

The following weighted $L^p$ bounds for elements of the Hilbert space $  D^{1}(\O)  $ can be found in \cite[Theorem B.1]{L2}.
\begin{lemma} \la{1leo}
   For   $m\in [2,\infty)$ and $\theta\in (1+m/2,\infty),$ there exists a positive constant $C$ such that we have for all $v\in  D^{1,2}(\O),$ \be\la{3h} \left(\int_{\O} \frac{|v|^m}{e+|x|^2}(\log (e+|x|^2))^{-\theta}dx  \right)^{1/m}\le C\|v\|_{L^2(B_1)}+C\|\na v\|_{L^2(\O) }.\ee
\end{lemma}

The following lemma was deduced in \cite{lx1}, we only state it here without proof.

\begin{lemma}\la{lemma2.6} For $\bar x$   as in \eqref{2.07},
suppose that $\n  \in L^\infty(\O)$ is a   function such that
\be \la{2.12}  0\le \n\le M_1, \quad M_2\le \int_{\b}\n dx ,\quad \n \bar x^\alpha \in L^1(\r^2),\ee
for $ N_*\ge 1 $ and positive constants $   M_1,M_2, $  and   $\al.$  Then, for $r\in [2,\infty),$ there exists a positive constant $C$ depending only on $  M_1, M_2, \alpha,   $ and $ r$  such that
 \be\la{z.1}\left(\int_{\r^2}\n |v |^r dx\right)^{1/r}  \le C  N_*^3  (1+\|\n\bar x^\al\|_{L^1(\r^2)})  \left(  \|\n^{1/2} v\|_{L^2(\r^2)} + \|\na  v \|_{L^2(\r^2)}\right) ,\ee for each $v\in \left.\left\{v\in D^1 (\O)\right|\n^{1/2}v\in L^2(\r^2) \right\}.$

\end{lemma}

Next, symbols $ \nabla^{\perp}\triangleq (-\p_2,\p_1),$ $\dot f\triangleq f_t+u\cdot\nabla f,$ denoting the material derivative of $f $. We state some elementary estimates which follow from (\ref{g1}) and the standard $L^p$-estimate  for the following elliptic system derived from the momentum equations in (\ref{a1}):
\be\la{h13}
\triangle F = \text{div}\left(\n\dot{u}-\text{div}(H\otimes H)\right),\quad\mu \triangle \o =
\nabla^\perp\cdot\left(\n\dot{u}-\text{div}(H\otimes H)\right) ,
\ee
where $F$ and $\omega$ are as in (\ref{hj1}).

\begin{lemma} \la{le3}
  Let $(\n,u,H)$ be a smooth solution of
   (\ref{a1}).
    Then for   $p\ge 2$ there exists a   positive
   constant $C$ depending only on $p,\mu,$ and $\lambda$ such that
\begin{eqnarray}
    &&\|{\nabla F}\|_{L^p(\O)} + \|{\nabla \o}\|_{L^p(\O)}
   \le C(\norm[L^p(\O)]{\n\dot{u}}+\norm[L^p(\O)]{|H||\na H|}),\label{h19}\\
&&\norm[L^p(\O)]{F} + \norm[L^p(\O)]{\o}
   \le C \left(\norm[L^2(\O)]{\n\dot{u}}+\norm[L^2(\O)]{|H||\na H|}\right)^{1-2/p }\nonumber\\
   &&\quad\cdot\left(\norm[L^2(\O)]{\nabla u} + \norm[L^2(\O)]{P }+\|H\|_{L^4}^2\right)^{2/p} ,
\la{h20}
\\
  &&\norm[L^p(\O)]{\nabla u} \le C \left(\norm[L^2(\O)]{\n\dot{u}}+\norm[L^2(\O)]{|H||\na H|}\right)^{1-2/p }\nonumber
   \\&&\quad\cdot\left(\norm[L^2(\O)]{\nabla u}
   + \norm[L^2(\O)]{P }+\|H\|_{L^4}^2\right)^{2/p}+ C\norm[L^p(\O)]{P}+C\||H|^2\|_{L^p}. \la{h18}
\end{eqnarray}
\end{lemma}
{\it Proof.} On one hand, the standard $L^p$-estimate for the elliptic system
(\ref{h13}) yields (\ref{h19}) directly, which, together with
(\ref{g1}) and (\ref{hj1}), gives (\ref{h20}).
On the other hand, since $-\Delta u=-\na {\rm div}u -\na^\perp\o,$ we have \bn\la{kq1}\na u=-\na(-\Delta)^{-1}\na {\rm
div}u-\na(-\Delta)^{-1}\na^\perp \o.\en Thus applying the standard $L^p$-estimate to \eqref{kq1} shows  \bnn \ba \|\na u\|_{L^p(\O)}&\le C(p) (\|{\rm
div}u\|_{L^p(\O)}+\|\o\|_{L^p(\O)})\\ &\le C (p)  \norm[L^p(\O)]{F} +C(p)  \norm[L^p(\O)]{\o} +
   C(p)  \norm[L^p(\O)]{P }+C\||H|^2\|_{L^p},\ea  \enn which,
along with (\ref{h20}), gives (\ref{h18}). Then, the proof of Lemma \ref{le3} is completed.

Next,  in order to get the
uniform (in time) upper bound of the density $\n,$ we need the following Zlotnik  inequality.
\begin{lemma}[\cite{zl1}]\la{le1}   Let the function $y$ satisfy
\bnn y'(t)= g(y)+b'(t) \mbox{  on  } [0,T] ,\quad y(0)=y^0, \enn
with $ g\in C(\r)$ and $y,b\in W^{1,1}(0,T).$ If $g(\infty)=-\infty$
and \be\la{a100} b(t_2) -b(t_1) \le N_0 +N_1(t_2-t_1)\ee for all
$0\le t_1<t_2\le T$
  with some $N_0\ge 0$ and $N_1\ge 0,$ then
\bnn y(t)\le \max\left\{y^0,\overline{\zeta} \right\}+N_0<\infty
\mbox{ on
 } [0,T],
\enn
  where $\overline{\zeta} $ is a constant such
that \be\la{a101} g(\zeta)\le -N_1 \quad\mbox{ for }\quad \zeta\ge \overline{\zeta}.\ee
\end{lemma}

Finally,    the following Beale-Kato-Majda type inequality,
which was proved in \cite{bkm} when $\div u\equiv 0,$   will be
used later to estimate $\|\nabla u\|_{L^\infty}$ and
$\|\nabla\n\|_{L^2\cap L^q} (q>2)$.
\begin{lemma}   \la{le9}  For $2<q<\infty,$ there is a
constant  $C(q)$ such that  the following estimate holds for all
$\na u\in L^2(\O)\cap D^{1,q}({\O }),$ \bnn \la{ww7}\ba \|\na
u\|_{L^\infty({\O })}&\le C\left(\|{\rm div}u\|_{L^\infty({\O })}+
\|\o\|_{L^\infty({\O })} \right)\log(e+\|\na^2
u\|_{L^q({\O })}) +C\|\na u\|_{L^2(\O)}+C. \ea\enn
\end{lemma}

\section{\la{se3} A priori estimates(I): lower order estimates}

In this section, we will establish some necessary a priori bounds
for smooth solutions to the Cauchy problem (\ref{a1})-(\ref{h2}) to extend the local strong   solution guaranteed by
Lemma \ref{th0}. Thus, let $T>0$ be a fixed time and $(\n,u,H)$ be
the smooth solution to (\ref{a1})-(\ref{h2})  on
${\r}^2 \times (0,T]$  with smooth initial
data $(\n_0,u_0,H_0)$ satisfying (\ref{co1}) and (\ref{h7}).

Set $\si(t)\triangleq\min\{1,t \}.$  Define
 \begin{eqnarray}
A_1(T)&\triangleq&\sup_{0\le t\le T}\sigma\left( \|\nabla u\|_{L^2}^2+\|\na H\|_{L^2}^2\right) + \int_0^{T} \sigma
 \left(\|\n^{1/2}\dot{u} \|_{L^2}^2 +\|\triangle H\|_{L^2}^2\right) dt,\la{As1}\\[2mm] A_2(T)&\triangleq&\sup_{0\le t\le T}\sigma^2\left(\|\n^{1/2}\dot{u} \|_{L^2}^2 +\||H||\na H| \|_{L^2}^2\right) \nonumber\\
  &&\quad+ \int_0^{T}\sigma^2\left(\|\na \dot{u} \|_{L^2}^2
  +\||\Delta H| |H|\|_{L^2}^2\right)dt,\la{As2}\\[2mm]
A_3(T)&\triangleq&\sup_{0\le t\le T}\sigma^{(3-2\beta)/4}(\|\na u\|_{L^2}^2+\|\na H\|_{L^2}^2).
\la{A3}
\end{eqnarray}

We have the following key a priori estimates on $(\n,u,H)$.
\begin{proposition}\la{pr1}  Under  the conditions of Theorem \ref{th1},
     there exists some  positive constant  $\ve$
    depending    on  $\mu ,  \lambda ,  \nu,  \ga ,  a ,  \on, \beta,$ $N_0,$ and $M$  such that if
       $(\n,u,H)$  is a smooth solution of
       (\ref{a1})-(\ref{h2})  on $\O \times (0,T] $
        satisfying
 \be\la{z1}
 \sup\limits_{
 \O \times [0,T]}\n\le 2\bar{\n},\quad
     A_1(T) + A_2(T) \le 2C_0^{1/2},\quad A_3(\sigma(T)) \leq 2C_0^{\delta_0},
   \ee
where $\delta_0=(2\beta-1)/(9\beta)$, the following estimates hold
        \be\la{z2}
 \sup\limits_{\O \times [0,T]}\n\le 7\bar{\n}/4, \quad
     A_1(T) + A_2(T) +\int_0^T\si \|P\|_{L^2}^2dt\le  C_0^{1/2}, \quad A_3(\si(T))\le C_0^{\delta_0},
  \ee
   provided $C_0\le \ve.$
\end{proposition}

The proof of Proposition \ref{pr1} will be postponed to the end of this section.

In the following, we will use the convention that $C$ denotes a generic positive constant depending  on $\mu$, $\lambda$, $\nu$, $\gamma$, $a$, $\bar{\n},$  $\beta,$ $N_0,$ and $M$, and  use $C(\al)$ to emphasize that $C$ depends on $\al.$


First, We will prove the  preliminary $L^2$ bounds for  $\nabla u$, $\na H$ and
$\n\dot{u}$.

\begin{lemma}\la{le21} Let $(\n,u,H)$ be a smooth solution of
 (\ref{a1})-(\ref{h2}) on $\O \times (0,T]. $ Then there exists a positive constant depending only on $\mu$,  $\lambda$,  $\nu$,  and $\gamma$ that
 \be \la{a16} \sup_{0\le t\le T}\int\left(\frac{1}{2}\n|u|^2+\frac{P}{\ga-1} +\frac{|H|^2}{2}\right)dx+\ia\int\left( \mu|\na
u|^2 +\nu|\na H|^2\right) dxdt\le C_0,\ee\be\la{h14}
  A_1(T) \le  C C_0 +C\sup\limits_{0\le t\le T}\|P\|_{L^2}^2 + C\int_0^{T}\sigma\int\left(|\nabla u|^3  +P|\nabla u|^2\right)dx dt,
  \ee
  \be\la{h15}
    A_2(T)
    \le   C C_0 +CA_1(T) +CA_1^2(T) + C\int_0^{T} \sigma^2 \left(\|\nabla u\|_{L^4}^4 +\|P\|_{L^4}^4\right) dt.
   \ee
\end{lemma}

{\it Proof.} First, the combination of standard energy inequality with \eqref{h3} gives \eqref{a16} directly.

Next, we will prove  \eqref{h14}.
Multiplying $(\ref{a1})_2 $ by
$  \dot{u}  $ and then integrating the resulting equality over
${\O } $ lead  to \be\la{m0} \ba    \int  \n|\dot{u} |^2dx      &
= - \int  \dot{u}  \cdot\nabla Pdx + \mu \int \triangle
 u\cdot \dot{u}  dx + (\mu+\lambda)
 \int \nabla\text{div}u\cdot \dot{u}  dx \\
 &\quad-\frac{1}{2}\int\dot u\cdot\na |H|^2 dx+\int H\cdot\na H\cdot\dot u dx\triangleq\sum_{i=1}^5R_i. \ea \ee

 Similar to the proof of   \cite[Lemma 3.2]{lx1}, we have
   \be\la{m1} \ba
\sum_{i=1}^3R_i\le &\left(\int  \text{div}uPdx -\frac{\mu }{2} \|\nabla u\|_{L^2}^2-\frac{\lambda+\mu}{2}
\|\text{div}u\|_{L^2}^2\right)_t\\& +C \int P|\na u |^2  dx+C\|\na u\|_{L^3}^3. \ea \ee

Using \eqref{a1}$_3$, we get
\be\label{lib1}\ba  
R_4&=\frac{1}{2}\int |H|^2  \div u_tdx+\frac{1}{2}\int |H|^2  \div (u\cdot \na u)dx\\
&=\left(\int \frac{|H|^2}{2} \div udx\right)_t+\frac{1}{2}\int u\cdot\na |H|^2 \div u dx +\frac{1}{2}\int |H|^2 \div(u\cdot \na u)dx\\
&\qquad-\int(H\cdot \na u+\nu \Delta H-H\div u)\cdot H\div udx\\
&=\left(\int \frac{|H|^2}{2} \div udx\right)_t-\frac{1}{2}\int |H|^2 (\div u)^2 dx+\frac{1}{2}\int |H|^2 \na u\cdot \na u dx\\
&\qquad-\int(H\cdot \na u+\nu \Delta H-H\div u)\cdot H\div u dx\\
&\leq \left(\int \frac{|H|^2}{2} \div udx\right)_t+C(\ve)\int |H|^2|\na u|^2dx+\ve \|\Delta H\|_{L^2}^2\\
&\leq \left(\int \frac{|H|^2}{2} \div udx\right)_t+C(\ve)\|\na u\|_{L^3}^2\|H\|_{L^2}^{4/3}\|\Delta H\|_{L^2}^{2/3}+\ve \|\Delta H\|_{L^2}^2\\
&\leq \left(\int \frac{|H|^2}{2} \div udx\right)_t+C(\ve)\|\na u\|_{L^3}^3 +2\ve \|\Delta H\|_{L^2}^2.
\ea\ee
Similarly, we have
\be\la{lib2}
\ba
R_5\le - \frac{\rm d}{{\rm d}t} \int H\cdot\na u\cdot H dx + C(\ve)\|\na u\|_{L^3}^3  +2\varepsilon\|\Delta H\|_{L^2}^2.
\ea
\ee

Putting \eqref{m1}-\eqref{lib2} into  \eqref{m0} yields
\be\la{lbq-jia20}\ba
& B'(t)
  +\int  \n|\dot{u} |^2 dx\le
  C \int P|\na u |^2  dx
 + C  (\ve)\|\nabla u\|_{L^3}^3 +4\varepsilon\|\Delta H\|_{L^2}^2,
 \ea\ee
 where \be\ba \la{nv1} B(t)\triangleq &\frac{\mu  }{2}\|\nabla
u\|_{L^2}^2 +\frac{ \lambda+\mu
}{2}\|\text{div}u\|_{L^2}^2-\int  \text{div}u P dx\\&-\frac{1}{2}\int\text{div}u|H|^2 dx+\int H\cdot\na u\cdot H dx.  \ea\ee

Next,
multiplying \eqref{a1}$_3$ by $ \triangle H$, and integrating by parts over $\mr^2$,  we have
\be\la{mm1}\ba
&\quad\frac{\rm d}{{\rm d}t}\int|\na H|^2 dx+2\nu\int|\triangle H|^2dx\\
&\quad\leq C\int |\na u||\na H|^2dx+ C\int |\na u|| H| |\Delta H|dx \\
&\quad\leq C\|\na u\|_{L^3} \|\na H\|_{L^2}^{4/3} \|\Delta H\|_{L^2}^{2/3}+C\|\na u\|_{L^3}\| H\|_{L^2}^{2/3} \|\Delta H\|_{L^2}^{4/3} \\
&\quad\leq C \|\na u\|_{L^3}^3+C \|\na H\|_{L^2}^4+\frac\nu2 \|\Delta H\|_{L^2}^{2}.
\ea\ee

Choosing $\tilde{C}$ suitably large such that
\be \la{n2'}\ba&\frac{\mu }{4}\|\nabla u\|_{L^2}^2+ \|\na H\|_{L^2}^2 -C \|P\|_{L^2}^2\\
  & \le B(t)+\tilde{C}\|\na H\|_{L^2}^2\le  C \|\nabla u\|_{L^2}^2+C\|\na H\|_{L^2}^2+ C \|P\|_{L^2}^2,\ea\ee
 adding  \eqref{mm1} multiplied by $\tilde{C}$  to \eqref{lbq-jia20}, and choosing $\ve$ suitably small lead to
\be\la{n1}\ba
&( B(t)+\tilde{C}\|\na H\|_{L^2}^2)' +\int \left( \n|\dot{u} |^2 +\nu\tilde{C}|\triangle H|^2\right)dx \\
& \le C \int P|\na u |^2  dx+ C  \|\nabla u\|_{L^3}^3+ C  \|\nabla H\|_{L^2}^4. \ea\ee
Hence,
we obtain (\ref{h14}) after integrating \eqref{n1} multiplied by $\si$
over $(0,T)  $ and using   \eqref{a16}  and \eqref{n2'}.

   Now, we will prove \eqref{h15}.  Operating $\pa/\pa_t+\div(u\cdot~)$ to $ (\ref{a1})_2^j $ and multiplying the resulting equation by $\dot{u}^j$, one gets by some simple calculations that
\be\label{lv3.40}\ba   \frac{1}{2}\left(\int \n |\dot{u}^j|^2dx\right)_t
&=\mu\int\dot{u}^j(\Delta u^j_t+\div(u\Delta u^j))dx\\&\quad+(\mu+\lambda)\int\dot{u}^j(\pa_t\pa_j(\div u)+\div(u\pa_j(\div u)))dx\\
&~~~-\int\dot{u}^j(\pa_jP_t+\div(u\pa_jP))dx\\&\quad-\frac{1}{2}\int\dot{u}^j(\pa_t\pa_j|H|^2+\div(u\pa_j|H|^2))dx\\
&~~~+\int\dot{u}^j(\pa_t(H\cdot\na H^j)+\div(u(H\cdot\na H^j)))dx\triangleq\sum^5_{i=1}I_i.
\ea\ee

 First, following the same arguments as in \cite{H3}, we  have
\be\label{lv3.41}\ba I_1+I_2+I_3&\leq -\frac{3\mu}{4}\|\na\dot{u}\|_{L^2}^2  +C\|\na u\|_{L^4}^4 +C\|P\|_{L^4}^4.
\ea\ee
Next, it follows from  \eqref{a1}$_3$ and \eqref{g1} that
\be\label{lv3.44}\ba I_4&= \int\pa_j\dot{u}^j H\cdot H_tdx +\frac{1}{2}\int\pa_i\dot{u}^j u^i\pa_jH^2 dx\\
&= \frac{1}{2}\int\pa_j\dot{u}^j  \pa_iu^iH^2dx -\frac{1}{2}\int\pa_i\dot{u}^j \pa_ju^iH^2 dx\\
&\quad+ \int\pa_j\dot{u}^j H\cdot (H\cdot \na u+\nu \Delta H-H\div u)dx\\
&\leq C\int|\na \dot{u}| |\na u||H|^2dx+C\int|\na \dot{u}| |\Delta H\cdot H| dx\\
&\leq \frac{\mu}{8}\int|\na \dot{u}|^2dx+C\int|\na u|^4dx+C\int |H|^8dx +C\int |\Delta H|^2|H|^2 dx.
\ea\ee
Similar to \eqref{lv3.44}, we estimate $I_5$ as follows 
\be\label{lv3.45}\ba I_5 \leq \frac{\mu}{8}\int|\na \dot{u}|^2dx+C\int|\na u|^4dx+C\int |H|^8dx +C\int |\Delta H|^2|H|^2 dx.\ea\ee
Putting \eqref{lv3.41}-\eqref{lv3.45} into \eqref{lv3.40}  yields
\be\label{lv3.46}\ba & \left(\int \n |\dot{u}|^2dx\right)_t+\mu\int|\na\dot{u}|^2dx\\
& \leq C \|\na u\|^4_{L^4}+C \|P\|^4_{L^4}+C \||H|^2\|^4_{L^4}+C \||\Delta H| |H|\|^2_{L^2}.
\ea\ee

Next,   we   estimate   the last term of \eqref{lv3.46}. For $ a_1,a_2\in\{-1,0, 1\},$ denote
 \be\label{amss1}\ba  \tilde{H}(a_1,a_2)=a_1H^1+a_2H^2,\quad\tilde{u}(a_1,a_2)=a_1u^1+a_2u^2.\ea\ee
It thus follows from \eqref{a1}$_3$ that
 \be\label{amss2}\ba  \tilde{H}_t-\nu\Delta \tilde{H}=H\cdot\na\tilde{u}-u\cdot \na\tilde{H}+\tilde{H}\div u .\ea\ee
 Integrating \eqref{amss2} multiplied by $4\nu^{-1}\tilde{H} \triangle |\tilde{H}|^2$  over $\mr^2$  leads to
\be\label{lv3.51}\ba & \nu^{-1}\left(  \|\na |\tilde{H}|^2\|^2_{L^2}\right)_t+ {2}\|\Delta |\tilde{H}|^2\|^2_{L^2}\\&=4\int  |\na \tilde{H}|^2 \Delta |\tilde{H}|^2dx-4\nu^{-1}\int H\cdot \na \tilde{u}\tilde{H}\Delta |\tilde{H}|^2dx\\
&\quad+4\nu^{-1}\int \div u |\tilde{H}|^2 \Delta |\tilde{H}|^2dx+ 2\nu^{-1}\int u\cdot \na |\tilde{H}|^2 \Delta |\tilde{H}|^2dx \\&
\le  C  \|\na u\|^4_{L^4}+C \|\na H\|^4_{L^4}+C  \||H|^2\|^4_{L^4}+ \|\Delta |\tilde{H}|^2\|^2_{L^2}  ,
\ea\ee
where we have used the following simple fact that
\bnn\ba
 \frac{1}{2}\int u\cdot \na |\tilde{H}|^2 \Delta |\tilde{H}|^2dx &=-\frac{1}{2}\int \na u\cdot \na |\tilde{H}|^2 \cdot\na |\tilde{H}|^2dx+\frac{1}{4}\int \div  u  |\na |\tilde{H}|^2|^2  dx\\
&   \leq C\|\na u\|^4_{L^4}+C\|\na H\|^4_{L^4}+C\||H|^2\|^4_{L^4}.
\ea\enn
Multiplying \eqref{lv3.51} by $\si^2$
and integrating the resultant inequality over $(0,T)$ lead to
 \be\label{nlv4'}\ba &\sup_{0\le t\le T}\si^2 \|\na |\tilde{H}|^2\|^2_{L^2}  +\int_0^T\si^2 \|\Delta |\tilde{H}|^2\|^2_{L^2}dt\\&\leq  CC_0+C\int_0^T\si^2\left(\|\na u\|^4_{L^4}+ \|\na H \|^4_{L^4}+ \||H|^2 \|^4_{L^4}\right)dt\\&\leq  CC_0+CA_1^2(T)+C\int_0^T\si^2 \|\na u\|^4_{L^4} dt+C\int_0^T\si^2 \||H|^2\|^4_{L^4} dt.
\ea\ee
To estimate the last term on the right hand side of \eqref{nlv4'}, we integrate  \eqref{a1}$_3$ multiplied by $|H|^2H$  over $\mr^2$  to obtain
\be\la{bug2}\ba
&\frac{1}{4}\frac{\rm d}{{\rm d}t}\||H|^2\|_{L^2}^2+\frac{\nu}{2}\|\na |H|^2\|_{L^2}^2+\nu\||H||\na H| \|_{L^2}^2 \\ & \leq C\|\na u\|_{L^2}\||H|^2\|_{L^4}^2 \leq  \frac{\nu}{4}\|\na |H|^2\|_{L^2}^2+C\||H|^2\|_{L^2}^2\|\na u\|_{L^2}^2,
\ea\ee
which together with   Gronwall's inequality and \eqref{a16} shows
\be\la{c1}\sup_{0\le t\le T}t\|H \|_{L^4}^4+\int_0^T t\||H||\na H| \|_{L^2}^2 dt\leq C\int_0^{T}\|H \|_{L^4}^4 dt\leq C C_0^{2}\ee
due to  \eqref{g1}. This combined with \eqref{g1}   gives  \be \la{cc1} \int_0^Tt^2\||H|^2\|^4_{L^4} dt\le C \int_0^Tt^2\|H\|^4_{L^4} \||H||\na H|\|_{L^2}^2dt\le CC_0^4.\ee

Noticing that
\be\label{lv3.57'1}\ba
 \||\na H||H| \|^2_{L^2}\le &  \|\na|\tilde{H}(1,0)|^2\|^2_{L^2}+ \|\na|\tilde{H}(0,1)|^2\|^2_{L^2} \\
 &+ \|\na|\tilde{H}(1,1)|^2\|^2_{L^2}+ \|\na|\tilde{H}(1,-1)|^2 \|^2_{L^2},\ea\ee and that \be\label{lv3.571}\ba
    \||\Delta H||H| \|^2_{L^2}\le  &  C \|\na H\|^4_{L^4}+  \|\Delta |\tilde{H}(1,0)|^2\|^2_{L^2}+ \|\Delta |\tilde{H}(0,1)|^2\|^2_{L^2}\\
& + \|\Delta |\tilde{H}(1,1)|^2\|^2_{L^2}+  \|\Delta |\tilde{H}(1,-1)|^2\|^2_{L^2},
\ea\ee we deduce from \eqref{nlv4'} and \eqref{cc1} that
 \be\label{nlv4}\ba &\sup_{0\le t\le T}\left( \si^2\||\na H||H| \|^2_{L^2}\right)+\int_0^T\si^2\left(\|\Delta |H|^2\|^2_{L^2}+\||\Delta H||H| \|^2_{L^2}\right)dt\\&\leq  CC_0+CA_1^2(T)+C\int_0^T\si^2 \|\na u\|^4_{L^4} dt.
\ea\ee

Finally,  multiplying  \eqref{lv3.46}  by $\sigma^2 $,  we obtain  (\ref{h15}) after using \eqref{nlv4}, \eqref{c1}, and  \eqref{cc1}.
The proof of Lemma \ref{le21} is    completed.

The next result shows that  pressure decays in time.
\begin{lemma}\la{le5} Let $(\n,u,H)$ be a smooth solution  of
   (\ref{a1})-(\ref{h2})     on $\O \times (0,T] $ satisfying (\ref{z1}). Then there exists a positive constant $C(\on)$ depending only  on $\mu,$  $\lambda,$  $\nu,$  $\gamma,$  and $\on$
 such that
  \be\la{h27}
  A_1(T)+A_2(T)+\int_0^T\si\|P\|_{L^2}^2dt\le C(\on) C_0 .
  \ee
   \end{lemma}

{\it Proof.} First, it follows from (\ref{h18}),   \eqref{z1},  \eqref{a16}, \eqref{c1}, and   \eqref{cc1}  that
\be\la{h99} \ba
  &  \int_0^{T}\sigma^2  \left(\|\na u\|_{L^4}^4 +\|P\|_{L^4}^4 \right) dt\\
& \le  C \int_0^{T}\si \left( \|\n  \dot u \|_{L^2}^2+\||H||\na H|\|_{L^2}^2\right)\left(\si\|\na u\|_{L^2}^2+\si\|P\|_{L^2}^2+\si\|H\|_{L^4}^4\right)dt\\
&\quad+C\int_0^T\si^2\||H|^2\|_{L^4}^4dt +C  \int_0^{T}\sigma^2  \|P\|_{L^4}^4  dt\\
& \le  C(\on)\left(A_1 (T)+C_0 \right)\int_0^{T}\left(\sigma  \|\n^{1/2}  \dot u \|_{L^2}^2+\si\||H||\na H|\|_{L^2}^2\right)dt\\
&\quad+CC_0^4 +C (\on) \int_0^{T}\sigma^2  \|P\|_{L^2}^2dt\\
& \le  C(\on)   C_0    +C (\on) \int_0^{T}\sigma^2  \|P\|_{L^2}^2dt.\ea \ee
To estimate the last term on the right-hand side of \eqref{h99}, we rewrite  $(\ref{a1})_2$ as
\be\la{u} P=(-\Delta )^{-1}\div(\n  \dot u )+(2\mu+\lambda)\div u+(-\Delta)^{-1}\div\div((H\otimes H)-\frac{1}{2}|H|^2),\ee
which together with \eqref{a16}, H\"older's  and Sobolev's inequalities yields  that
\bnn \la{u'} \ba \int P^2dx
\le&C \|(-\Delta)^{-1} \div(\n  \dot u )\|_{L^{4\ga}} \|P\|_{L^{4\ga/(4\ga-1)}} +C\|\na u\|_{L^2}\|P\|_{L^2}+C\||H|^2\|_{L^2}\|P\|_{L^2}  \\
\le& C\|  \n  \dot u \|_{L^{4\ga/(2\ga+1)}}\|\n\|_{L^1}^{1/2}\|\n\|_{L^{2\ga}}^{\ga-1/2}+C\|\na u\|_{L^2} \|P\|_{L^2}+C\|H\|_{L^2}\|\na H\|_{L^2}\|P\|_{L^2}\\
 \le& C \|\n^{1/2}\|_{L^{4\ga}} \| \n^{1/2}  \dot u  \|_{L^2}\|\n\|_{L^1}^{1/2}\|\n\|_{L^{2\ga}}^{\ga-1/2} +C\|\na u\|_{L^2} \|P\|_{L^2}+C\|\na H\|_{L^2}\|P\|_{L^2}\\
 \le& C\|P\|_{L^2}  \| \n^{1/2}  \dot u  \|_{L^2} +C\|\na u\|_{L^2} \|P\|_{L^2}+C\|\na H\|_{L^2}\|P\|_{L^2},\ea \enn
where in the last inequality, one has used
  \be \la{mr3}\int\n dx =\int \n_0dx =1,\ee
 due to   the mass conservation equation  $\eqref{a1}_1.$
Thus, we arrive at
\be \la{new1}\ba \|P\|_{L^2}  \le   C  \| \n^{1/2}  \dot u \|_{L^2} + { C}\|\na u\|_{L^2}+C\|\na H\|_{L^2},\ea \ee
which, along with  \eqref{h14},   \eqref{h15}, \eqref{h99}, \eqref{a16}, and \eqref{z1} gives   \be\la{h28}\ba
A_1(T)+A_2(T)\le& C(\on)C_0 +C(\on)\int_0^{T}  \sigma\|\nabla u\|_{L^3}^3 dt. \ea\ee

 Finally, on the one hand,  one deduces from (\ref{h18}),  \eqref{g1}, (\ref{a16}),   (\ref{c1}),  and (\ref{z1})  that \be\la{h34} \ba
 & \int_0^{\si(T)} \sigma\|\nabla u\|_{L^3}^3 dt\\
&  \le C \int_0^{\si(T)} \sigma\left(\|\n^{1/2} \dot u \|_{L^2}+\||H||\na H|\|_{L^2} \right) \left(\|\na u\|_{L^2}^2+\|P\|_{L^2}^2 +\|H\|_{L^4}^4\right)  dt \\
&\quad+C\int_0^{\si(T)}\si\|P\|_{L^3}^3 dt +C\int_0^{\si(T)}\si \|\na |H|^2\|_{L^2}^2\||H|^2\|_{L^2}^2dt\\
&  \le C A_2^{1/2}(\si(T))\int_0^{\si(T)}  \left(\|\na u\|_{L^2}^2+\|P\|_{L^2}^2+\|\na H\|_{L^2}^2 \right)dt+C\int_0^{\si(T)}\|P\|_{L^3}^3 dt\\
& \le  C(\on)C_0.\ea \ee
%
 On the other hand, H\"older's inequality, \eqref{h99},  \eqref{z1}, and \eqref{new1} imply \be\la{h33} \ba
 \int_{\si(T)}^{T}  \|\nabla u\|_{L^3}^3  dt
&\le
\de\int_{\si(T)}^{T}    \|\nabla u\|_{L^4}^4 dt + C(\de) \int_{\si(T)}^{T} \|\nabla u\|_{L^2}^2 dt \\& \le  \de C(\on) A_1(T)+  C(\de)C(\on) C_0.
\ea \ee
Putting   \eqref{h34}  and \eqref{h33} into \eqref{h28} and choosing $\de$ suitably small lead to \bnn A_1(T)+A_2(T)\le C(\on)C_0,\enn  which together with \eqref{new1} and \eqref{a16} gives \eqref{h27} and finishes  the proof of Lemma \ref{le5}.


Next, we derive the rates of decay for $\na u, \na H, H$ and $P$, which are  essential to obtain the uniform  (in time) upper bound of the density for large time.
\begin{lemma}\la{ly1}
    For $p\in [2,\infty),$  there exists a positive constant $C(p,\on)$
  depending only  on $p, \mu, \nu,$  $\lambda,$  $\gamma,$  and $\on$ such that,  if  $(\n,u,H)$ is a smooth solution  of
   \eqref{a1}--\eqref{h2}     on $\O \times (0,T] $ satisfying \eqref{z1}, then
\be \la{ly8} \ba &\sup\limits_{\si(T)\le t\le T}\left[t^{p-1}\left(\|\na u\|_{L^p}^p + \|P\|^p_{L^p}\right)
+t\left(\|\na H\|_{L^2}^2 + \||H|^2\|^2_{L^2}\right)\right]\\
&\quad+\sup\limits_{\si(T)\le t\le T}\left[t^2\left(\|\n^{1/2}\dot u\|_{L^2}^2 +\||\na H||H|\|_{L^2}^2\right) \right] \le C(\on)C_0.\ea\ee
\end{lemma}

{\it Proof.} First, the pressure  $P$ satisfies \be
\la{a95}P_t+u\cdot\nabla P+\ga P{\rm div}u=0 .\ee For $p\ge 2,$ multiplying \eqref{a95} by $p P^{p-1}$ and integrating
the resulting equality over ${\O } ,$ one gets after using ${\rm
div}u=\frac{1}{2\mu+\lambda}\left(F+P+\frac{1}{2}|H|^2\right)$ that \be\la{vv1}\ba&
\left(\| P\|_{L^p}^p \right)_t+ \frac{p\ga-1}{2\mu+\lambda}\|P\|_{L^{p+1}}^{p+1}
\\ &=- \frac{p\ga-1}{2\mu+\lambda}\int
P^p\left(F+\frac{1}{2}|H|^2\right)dx \\&
\le \frac{p\ga-1}{2(2\mu+\lambda)}\|P\|_{L^{p+1}}^{p+1}+C(p)  \|F\|_{L^{p+1}}^{p+1}+C(p)  \||H|^2\|_{L^{p+1}}^{p+1},\ea\ee
which together with \eqref{h20} and \eqref{g1} gives
\be\la{a96}\ba
 &\frac{2(2\mu+\lambda)}{p\ga-1}\left(\| P\|_{L^p}^p \right)_t+ \|P\|_{L^{p+1}}^{p+1}
   \le C(p) \|F\|_{L^{p+1}}^{p+1}+C(p)  \||H|^2\|_{L^{p+1}}^{p+1} \\
    &\le C(p) \left(\|\na u\|_{L^2}^2+\|P\|_{L^2}^2+\|H\|_{L^4}^4\right)\left(\|\n \dot u\|_{L^2}+\||H| |\na H|\|_{L^2}\right)^{p-1}.\ea\ee
In particular, choosing $p=2$  in \eqref{a96} shows
\be\ba
\la{a9z6}
&\left(\| P\|_{L^2}^2 \right)_t+ \frac{2\ga-1}{2(2\mu+\lambda)}\|P\|_{L^3}^3  \\
  &\leq \de \|\n^{1/2}\dot u\|_{L^2}^2+\de\||H||\na H|\|_{L^2}^2+C(\de)\left(\|\na u\|_{L^2}^4+\|\na H\|_{L^2}^4+\|P\|_{L^2}^4\right),\ea\ee
where in the last inequality we have used \eqref{g1} and \eqref{a16}.

Next, it follows from \eqref{n1}, \eqref{h18},  \eqref{g1}, and \eqref{a16} that
\be   \la{ly11}\ba
 &( B(t)+\tilde{C}\|\na H\|_{L^2}^2)' +\int \left( \n|\dot{u} |^2 +\nu|\triangle H|^2\right)dx\\
 & \le C\| P\|_{L^3}^3+ C  \|\nabla u\|_{L^3}^3+C\|\na H\|_{L^2}^4 \\
& \le C_1  \|P\|_{L^3}^3+ C \left(\|\n\dot u\|_{L^2}+\||H||\na H|\|_{L^2}\right)\left(\|\na u\|_{L^2}^2+\|P\|_{L^2}^2+\|H\|_{L^4}^4\right)\\
&\quad+C\||H|^2\|_{L^3}^3 +C\|\na H\|_{L^2}^4\\
& \le C_1 \|P\|_{L^3}^3+ \frac14\|\n^{1/2}\dot u\|_{L^2}^2+C\||H||\na H|\|_{L^2}^2 \\&\quad+C(\on )\left(\|\na u\|_{L^2}^4+\|\na H\|_{L^2}^4+\|P\|_{L^2}^4\right).
\ea\ee

 Choosing $ C_2\ge 2+  2(2\mu+\lm)(C_1+1)/(2\ga-1)$
suitably large such that
\be\la{ly12}\ba &\frac{\mu}{4}\|\na u\|_{L^2}^2+ \|\na H\|_{L^2}^2+\|P\|_{L^2}^2 \\&\le B(t)+\tilde{C}\|\na H\|_{L^2}^2+C_2\|P\|_{L^2}^2 \le C\|\na u\|_{L^2}^2+C\|\na H\|_{L^2}^2+C\|P\|_{L^2}^2,\ea\ee
and adding \eqref{a9z6} multiplied by $C_2$  to \eqref{ly11}, we obtain after choosing $\de$ suitably small that
\be \la{h29}\ba &2\left( B (t)+\tilde{C}\| \na H \|_{L^2}^2+C_2 \| P\|_{L^2}^2 \right)' +\int  \left(\n|\dot{u} |^2+\nu|\triangle H|^2 + P^3  \right)dx \\
& \le
   C \| P\|_{L^2}^4 + C  \|\na u\|_{L^2}^4 +C\|\na H\|_{L^2}^4+C\||H||\na H|\|_{L^2}^2  ,\ea\ee
   which  multiplied  by $t,$  together with Gronwall's inequality,   \eqref{ly12},    \eqref{h27},    \eqref{a16},  \eqref{c1},  and  \eqref{z1}    yields
 \be \la{h31} \ba & \sup\limits_{
\si(T)\le t\le T}t\left( \| \na u\|_{L^2}^2+\|\na H\|_{L^2}^2+ \| P\|_{L^2}^2+\| H \|_{L^4}^4
 \right) \\& + \int_{\si(T)}^Tt \int  \left(\n|\dot{u} |^2+|\triangle H|^2 +|H|^2|\na H|^2+ P^3  \right)dxdt\le C(\on)C_0. \ea\ee


Next,
multiplying \eqref{lv3.46} by $t^2 $ shows
\be\la{zo1} \ba \left(t^2
  \int\n|\dot{u}|^2dx \right)_t + { \mu}t^2 \int
 |\nabla\dot{u}|^2dx  & \le 2t\int\n|\dot{u}|^2dx+Ct^2\|\na u\|_{L^4}^4+Ct^2\|P\|_{L^4}^4 \\
 &\quad+Ct^2\||H|^2\|_{L^4}^4+C_3t^2\||\Delta H||H|\|_{L^2}^2.\ea\ee

Adding \eqref{lv3.51}  multiplied by $(4C_3+1)t^2 $    to \eqref{zo1} yields that for $t\in [\si(T),T]$
\be\la{lbq-jia25} \ba & \left(t^2
  \int\n|\dot{u}|^2dx+(4C_3+1) \nu^{-1}t^2\|\na |\tilde{H}|^2\|^2_{L^2}\right)_t\\
  &\quad+ {\mu}t^2\|\nabla\dot{u}\|^2_{L^2}+(4C_3+1)t^2\|\Delta |\tilde{H}|^2\|^2_{L^2} \\
 & \leq C_3t^2\||\Delta H||H|\|_{L^2}^2+Ct\left(\||\na H||H| \|^2_{L^2}+\|\n^{1/2}\dot{u}\|^2_{L^2}\right)
 \\
  &\quad+Ct^2\left(\|\na u\|^4_{L^4}+ \|\na H \|^4_{L^4}+ \||H|^2\|^4_{L^4}\right)+ Ct^2\|P \|^4_{L^4}\\
   & \leq C_3t^2\||\Delta H||H|\|_{L^2}^2+ Ct\left(\||\na H||H| \|^2_{L^2}+\|\n^{1/2}\dot{u}\|^2_{L^2}\right)\\
  &\quad +Ct^2\left(\|\n^{1/2}\dot{u}\|^2_{L^2}+\||\na H||H| \|^2_{L^2}\right)\left(\|\na u\|^2_{L^2}+ \|P \|^2_{L^2}+ \|H\|^4_{L^4}\right)\\
  &\quad+ Ct^2\|\na H \|^2_{L^2}\|\na^2 H \|^2_{L^2}  +
 C t^2\|P \|^4_{L^4}\\
   & \leq C_3t^2\||\Delta H||H|\|_{L^2}^2+Ct\left(\|\n^{1/2}\dot{u}\|^2_{L^2}+\||\na H||H| \|^2_{L^2}+\|\na^2 H \|^2_{L^2}\right) + C_4t^2\|P \|^4_{L^4}
 \ea\ee
where in the second inequality we have used \eqref{h18},    \eqref{g1}, and  \eqref{c1}.

Choosing $p=3$ in \eqref{a96}, we add  \eqref{a96}  multiplied by $(C_4+1)t^2$ to \eqref{lbq-jia25} and obtain that for $t\in [\si(T),T]$
\be\la{lbq-jia25'} \ba & \left(t^2
  \int\n|\dot{u}|^2dx+ (4C_3+1)\nu^{-1}t^2\|\na |\tilde{H}|^2\|^2_{L^2}+ \frac{2(2\mu+\lambda)( C_4 +1)}{3\ga-1}t^2\| P\|_{L^3}^3\right)_t\\
  &\quad+ {\mu}t^2\|\nabla\dot{u}\|^2_{L^2}+(4C_3+1)t^2\|\Delta |\tilde{H}|^2\|^2_{L^2}+t^2\|P\|^4_{L^4}\\
   & \leq C_3t^2\||\Delta H||H|\|_{L^2}^2+Ct\left(\|\n^{1/2}\dot{u}\|^2_{L^2}+\||\na H||H| \|^2_{L^2}+\|\na^2 H \|^2_{L^2}+\|P \|^3_{L^3}\right) \ea\ee
which combined with \eqref{h31},  \eqref{lv3.57'1}, \eqref{lv3.571}, \eqref{nlv4},  and   \eqref{z1} yields
\be \la{ly9} \ba  & \sup\limits_{\si(T)\le t\le T}t^2\int  \left(\n|\dot{u} |^2+ |\na H|^2|H|^2+ P^3  \right)dx \\
&\quad+  \int_{\si(T)}^Tt^2 \left(\|\nabla\dot{u}\|^2_{L^2}+\|\Delta |H|^2\|^2_{L^2}+\||\Delta H||H| \|^2_{L^2}+\|P\|_{L^4}^4\right)dt \le C(\on)C_0. \ea\ee

Finally, we claim that for   $m=1,2,\cdots,
$ \be \la{ly10}\sup\limits_{\si(T)\le t\le   T}t^{m}\|P\|_{L^{m+1}}^{m+1}+\int_{\si(T)}^Tt^m\|P\|_{L^{m+2}}^{m+2}dt\le C(m,\on)C_0,\ee which together with \eqref{h18}, \eqref{h31}, and \eqref{ly9} gives \eqref{ly8}.   We shall prove \eqref{ly10} by induction. In fact, for $m=1,$ \eqref{ly10} holds  due to \eqref{h31}. Assume that \eqref{ly10} holds for $m=n,$ that is, \be \la{ly16}\sup\limits_{\si(T)\le t\le T}t^{n}\|P\|_{L^{n+1}}^{n+1}+\int_{\si(T)}^Tt^n\|P\|_{L^{n+2}}^{n+2}dt\le C(n,\on)C_0.\ee
Multiplying \eqref{a96} where  $p=n+2$  by $t^{n+1} ,$ one obtains after  using \eqref{ly9}
\be\la{zo8}\ba &
  \frac{2(2\mu+\lambda)}{(n+2)\ga-1}\left(t^{n+1}\| P\|_{L^{n+2}}^{n+2} \right)_t+ t^{n+1}\|P\|_{L^{n+3}}^{n+3}\\
  &\le C (n,\on) t^{n}\| P\|_{L^{ n+2} }^{n+2} + Ct^{n+1}\||H|^2\|_{L^{n+3}}^{n+3} \\
&\quad+C(n,\on)\left(t\|\n^{1/2}\dot{u}\|_{L^2}+t\||\na H||H| \|_{L^2}\right)^{n+1}\left(\|\na u\|^2_{L^2}+ \|P \|^2_{L^2}+\||H|^2\|^2_{L^2}\right)\\
&\le C (n,\on) t^{n}\| P\|_{L^{ n+2} }^{n+2}+C(n,\on)C_0\left(\|\na u\|^2_{L^2}+ \|P \|^2_{L^2}+\|\na H\|^2_{L^2}\right).
\ea\ee
Integrating \eqref{zo8} over $[\si(T),T]$ together with \eqref{ly16} and  \eqref{h27} shows that \eqref{ly10} holds for $m=n+1.$   By induction, we obtain \eqref{ly10} and finish the proof of Lemma \ref{ly1}.

 Next, Lemma \ref{lemma2.6} combined with the following Lemma \ref{lemma2.7} which has been proved in \cite[Lemma 3.5]{lx1}, will be useful to estimate the $L^p$-norm of $\n \dot u$   and obtain
     the uniform  (in time)
upper bound of the density for large time.

\begin{lemma}\la{lemma2.7} Let $(\n,u,H)$ be a smooth solution of (\ref{a1})-(\ref{h2}) on $\O \times (0,T] $ satisfying the assumptions in Theorem \ref{th1} and (\ref{z1}). Then for any $\al>0,$ there
exists a    positive constant  $N_1$   depending only on   $\al,$     $N_0,$ and $M$ such
that for all $t\in (0,T],$
    \be\la{uq2}   \int_{B_{N_1(1+t)\log^\al(1+t)}}\n(x,t) dx \ge  \frac14. \ee

\end{lemma}

Next,    to  obtain      the
upper bound of the density for small time, we still need the following lemmas about the short-time estimates of $H$ and $u$.

\begin{lemma}\la{chuzhi} Let $(\n,u,H)$ be a smooth solution of (\ref{a1})-(\ref{h2})   on $\O \times (0,T] $ satisfying (\ref{z1}) and the assumptions in Theorem \ref{th1}. Then there
exists a positive constant  $C$   depending only on $\mu ,  \lambda ,  \nu, \ga ,  a ,  \on, \beta,$ $N_0,$ and $M$  such
that
   \be\la{lbqnew-czh}  \sup_{0\le t\le \si(T)}t^{1-\beta}  \|\na
H\|_{L^2}^2 +\int_0^{\si(T)}t^{1-\beta}(\|\na
H\|_{L^2}^2+\|\na^2
H\|_{L^2}^2)dt\le C(M). \ee
\end{lemma}

{\it Proof.}
First, integrating   \eqref{a1}$_3$ multiplied  by $H$   over $\mathbb{R}^2$ gives
\be\la{czh1}\ba
\frac{1}{2}\frac{\rm d}{{\rm d}t}\|H\|_{L^2}^2 +\nu\|\na H\|_{L^2}^2
 \le \frac{\nu}{2}\|\na H\|_{L^2}^2+\|\na u\|_{L^2}^2\|H\|_{L^2}^2
\ea\ee
which together with \eqref{a16} and Gronwall's inequality implies that
\be\la{czh2}\ba
 \sup_{0\le t\le T}\|H\|_{L^2}^2+\int_0^{T}\|\na
H\|_{L^2}^2dt\le C\|H_0\|_{L^2}^2.
\ea\ee

 Next, it follows from  \eqref{mm1} and the Gargliardo-Nirenberg inequality that
\be\la{lbq-jia31}\ba&
\frac{\rm d}{{\rm d}t}\int|\na H|^2 dx+ 2\nu\int|\triangle H|^2dx\\&\le C\int |\na u||\na H|^2dx+ C\int |\na u|| H| |\Delta H|dx\\
& \leq  C\|\na u\|_{L^2} \|\na H\|_{L^2}\|\na^2H\|_{L^2}+C\|\na u\|_{L^2} \|H\|_{L^2}^{1/2} \|\Delta H\|_{L^2}^{3/2} \\
 &\leq C \|\na u\|_{L^2}^2\|\na H\|_{L^2}^2+C \|\na u\|_{L^2}^{4}\|H\|_{L^2}^{2}+
\nu\|\Delta H\|_{L^2}^2,
\ea\ee
which together with \eqref{czh1} gives
\be\la{czh3}\ba
\frac{\rm d}{{\rm d}t} \|H\|_{H^1}^2 + \nu\|\na
H\|_{H^1}^2
\leq  C\left(\|\na u\|_{L^2}^2+\|\na u\|_{L^2}^{4}\right) \|H\|_{H^1}^2 .
\ea\ee
Noticing that   \eqref{z1} gives
  \be\la{chu2}\ba
&  \int_0^{\si(T)}(\|\na u \|_{L^2}^4+\|\na H\|_{L^2}^4)dt
\\
&\leq \sup_{  0\le t\le \si(T)}\left(\sigma^{(3-2\beta)/4}(\|\na u \|_{L^2}^2+\|\na H\|_{L^2}^2)\right)^2\int_0^{\si(T)}\sigma^{(2\beta-3)/2}dt  \\
&\leq CC_0^{2\delta_0}
\ea
\ee
due to $\beta\in(1/2,1],$  we obtain from \eqref{czh3}   by using   Gronwall's inequality   that
   \be\la{czh4}\ba  \sup_{0\le t\le\si( T)} \|H\|_{H^1}^2  +\int_0^{\si(T)} \|\na
H\|_{H^1}^2 dt\le C\|H_0\|_{H^1}^2. \ea\ee
Then, multiplying \eqref{czh3} by $t$ and integrating it over $(0,\si(T))$, by \eqref{czh2},   \eqref{a16} and \eqref{chu2}, it holds that
   \be\la{czh5}\ba  &\sup_{0\le t\le \si(T)}t \|H\|_{H^1}^2 +\int_0^{\si(T)}t \|\na H\|_{H^1}^2 dt\le C\|H_0\|_{L^2}^2. \ea\ee

Since the solution operator $H_{0}\mapsto H(\cdot,t)$ is linear,
by the standard Stein-Weiss interpolation argument (see \cite{bl}),
one can obtain \eqref{lbqnew-czh}  directly from (\ref{czh4})   and (\ref{czh5}).
Thus, we finish the proof of Lemma \ref{chuzhi}.

\begin{lemma}\la{zc1} Let $(\n,u,H)$ be a smooth solution of (\ref{a1})-(\ref{h2})   on $\O \times (0,T] $ satisfying (\ref{z1}) and the assumptions in Theorem \ref{th1}. Then there
exists a positive constant  $C$   depending only on $\mu ,  \lambda , \nu,  \ga ,  a ,  \on, \beta,$ $N_0,$ and $M$  such
that
   \be\la{uv1}  \sup_{0\le t\le \si(T)}t^{1-\beta} \|\na
u\|_{L^2}^2+\int_0^{\si(T)}t^{1-\beta}\int\n|\dot u|^2dxdt\le
C(\on,M). \ee
\end{lemma}

{\it Proof.}
As in \cite{hof2002},  for  the linear differential operator $L$
defined by \be\ba(Lw)^j&\triangleq \n w^j_{t}+ \n u\cdot\na w^j
-(\mu\Delta w^j+(\mu+\lambda) \pa_j\div w )\no&=\n\dot w^j
-(\mu\Delta w^j+(\mu+\lambda) \pa_j\div w ),\quad j=1,2, \ea\ee let $w_1$, $w_2$ and $w_3$   be
the solution to: \be \la{sas2}Lw_1 =0 ,\quad w_1(x,0)=w_{10}(x), \ee
\be\la{sas3} Lw_2 = -\nabla P(\n),\quad w_2(x,0)=0,\ee
and
\be\la{sasw3}
Lw_3=\frac{1}{2}\left(H\cdot\na B+B\cdot \na H\right)-\frac{1}{2}\left(\na B\cdot H+\na H\cdot B\right),\quad w_3(x,0)=0,
\ee
respectively, where $B=(B^1,B^2)$ is the solution of the following initial problem \be\la{a1'}\begin{cases}
 B_t-\na \times(u\times B)=-\na \times(\nu \na \times B),\quad \div B=0,\\ B(x,0)=B_0(x),\end{cases}\ee
with   initial data $B_0(x)$.

 First,
similar to \cite[Lemama 4.1]{chl}, multiplying  \eqref{a1}$_3$   by $H\bar{x}^a$ and integrating by parts yields
\be \la{lv4.1}  \ba
&\frac{1}{2}\left(\int |H|^2\bar{x}^adx\right)_t+\nu \int |\na H|^2\bar{x}^adx=\frac{\nu}{2}\int |H|^2\Delta\bar{x}^adx\\
&\quad+\int H\cdot\na u\cdot H\bar{x}^adx-\frac{1}{2}\int \div u|H|^2\bar{x}^adx+\frac{1}{2}\int |H|^2u\cdot\na\bar{x}^adx\triangleq\sum^4_{i=1}\bar{I}_i.
\ea\ee

Direct calculations yields that
\be \la{lv4.4}
 |\bar{I}_1| \leq C\int |H|^2 \bar{x}^a \bar{x}^{-2}\log^{4}(e+|x|^2) dx\leq C\int |H|^2 \bar{x}^a dx,\ee and that
 \be\la{ljo2}\ba
 |\bar{I}_2|+|\bar{I}_3|&\leq C\int |\na u||H|^2 \bar{x}^a dx\leq C \|\na u\|_{L^2}\| H   \bar{x}^{a/2}\|_{L^4}^2\\
 &\leq C \|\na u\|_{L^2}\| H   \bar{x}^{a/2}\|_{L^2} (\| \na H \bar{x}^{a/2}\|_{L^2}+\| H \na\bar{x}^{a/2}\|_{L^2})\\
 &\leq C(\|\na u\|_{L^2}^2+1)\| H \bar{x}^{a/2}\|_{L^2}^2 +\frac\nu4\| \na H \bar{x}^{a/2}\|_{L^2}^2.\ea\ee

 It follows from  \eqref{3h}, \eqref{uq2}, and the Poincar\'e-type inequality     \cite[Lemma 3.2]{F2} that for  $s>2, \delta\in(0,1], $ and $t\in [0,\si(T)],$\be \la{ljo1}\|u\bar x^{-1}\|_{L^2}+ \|u\bar x^{-\delta}\|_{L^{s/\delta}}\le C(s,\delta,\on,M)\|\n^{1/2}u\|_{L^2}+C(s,\delta,\on,M)\|\na u\|_{L^2},\ee which yields that
 \be\la{ljo3}\ba
|\bar{I}_4|&\leq    C\| H \bar{x}^{a/2}\|_{L^4}\| H \bar{x}^{a/2}\|_{L^2} \|u\bar{x}^{-3/4}\|_{L^{4}}\\
&\leq   C\| H \bar{x}^{a/2}\|_{L^4}^2+C\|H \bar{x}^{a/2}\|_{L^2}^2 \|u\bar{x}^{-3/4}\|_{L^{4}}^2 \\
&\leq C(\on,M)\left(1+ \|\na u\|_{L^2}^2\right)\|H\bar{x}^{a/2}\|_{L^2}^2+\frac\nu4\| \na H \bar{x}^{a/2}\|_{L^2}^2.
\ea\ee
 Putting \eqref{lv4.4}, \eqref{ljo2}, and  \eqref{ljo3} into \eqref{lv4.1}, after using Gronwall's inequality and \eqref{a16}, we have
\be
\ba\la{lbqnew-gj10}
\sup_{0\le t\le \sigma(T)} \int \bar{x}^{a}|H|^2 dx +\int_{0}^{\sigma(T)} \int  \bar{x}^{a}|\na H|^2 dxdt\le C(\|H_0\bar{x}^{a/2}\|_{L^2}^2)\le C(M),
\ea
\ee

Next,  set $$\tau\triangleq \min\left\{\frac{\mu^{1/2}}{ 2(1+2\mu+\lm)^{1/2} },\,\frac{ \beta}{1-\beta}\right\}\in (0,1/2]  .$$ If $\beta\in (\frac{1}{2},1),$   Sobolev's inequality implies
\be \la{zz.1}\ba\int \n_0|u_0|^{2+\tau}dx &\le  \int \n_0|u_0|^{2 }dx +  \int \n_0 |u_0|^{2/(1-\beta)}dx\\ &\le C(\on)+C(\on)\|u_0\|_{\dot H^\beta}^{2/(1-\beta)}\le C(\on,M).\ea\ee For the case that $\beta=1,$ one obtains from \eqref{z.1} that
\be  \la{zz.2}\ba\int \n_0|u_0|^{2+\tau}dx &\le C(\on)\left( \int \n_0|u_0|^{2 }dx +  \int   |\na u_0|^2dx\right)^{(2+\tau)/2} \le C(\on,M).\ea\ee

Then, multiplying $(\ref{a1})_2$ by $(2+\tau)|u|^\tau u$  and   integrating the resulting equation over $ \O$ lead  to
\be\la{lbq-jia1}\ba & \frac{d}{dt}\int \n |u|^{2+\tau}dx+ (2+\tau) \int|u|^\tau \left(\mu |\na u|^2+(\mu+\lm) (\div u)^2\right) dx \\& \le  (2+\tau)\tau  \int (\mu+\lm)|\div u||u|^\tau |\na u|dx +C  \int \n^\ga |u|^\tau |\na u|dx+C\int|H|^2|u|^\tau|\na u|dx \\& \le \frac{2+\tau}{2}\int (\mu+\lm)(\div u)^2|u|^\tau  dx +\frac{(2+\tau)\mu}{4} \int  |u|^\tau |\na u|^2 dx \\&\quad + C \int \n |u|^{2+\tau} dx+C\int \n^{(2+\tau)\ga-\tau/2}dx+C\int|H|^4 |u|^\tau dx.\ea\ee
It follows from \eqref{lbqnew-gj10}, \eqref{chu2}, and \eqref{ljo1}  that
\be\la{lbq-jia2}\ba&\int_0^{\si(T)} \int|H|^4 |u|^\tau dxdt\\&\leq C\int_0^{\si(T)}\||H|^{4-2\tau}\|_{L^{2/(2-3\tau)}}\|(|H|\bar x^{a/2})^{2\tau}\|_{L^{1/\tau}}\||u|^\tau\bar x^{-a\tau}\|_{L^{2/\tau}}dt\\
&\leq  C(\on,M)\int_0^{\si(T)} \|H\|^{4-2\tau}_{H^1}\left(\|\rho^\frac{1}{2}u\|_{L^2}+\|\na u\|_{L^2} \right)^{\tau}dt\\&\le C(\on,M)\int_0^{\si(T)}\left(1+\|\na  H \|_{L^2}^4+ \|\na u\|_{L^2}^2\right)dt\le C(\on,M).
\ea\ee
     Then applying Gronwall's inequality to \eqref{lbq-jia1}, together with \eqref{zz.1}, \eqref{zz.2},  and \eqref{lbq-jia2} yields that
   \be\la{jan2}
  \sup_{0\le t\le  \si(T) }\int \n |u|^{2+\tau}dx\le C(\on,M ).
  \ee

With \eqref{jan2} at hand,    similar to the proof of   \cite[Lemma 3.6]{lx1}, we can obtain that
\be \label{uu4} \sup_{0\le t\le
\si(T)}t^{1-\beta}\|\nabla
w_1\|_{L^2}^2+\int_0^{\si(T)}t^{1-\beta}\int\n|\dot
{w_1}|^2dxdt\leq C(\on,M)\| w_{10}\|_{\dot H^\beta}^2 ,\ee  and that
\be  \label{uu3} \sup_{0\le t\le \si(T)}\|\nabla
w_2\|_{L^2}^2+\int_0^{\si(T)}\int\n|\dot {w_2}|^2dxdt\leq
C(\on,M).\ee

Finally, it remains to  estimate $w_3.$

 First, multiplying \eqref{sasw3} by $\dot{w}_3=w_{3t}+u\cdot\na w_3$  and integrating by parts over $\mathbb{R}^2$ lead to \be\la{czw6}\ba
&\frac{\mu}{2}\frac{\rm d}{{\rm d}t}\int|\na w_3|^2dx+\int\rho|\dot w_3|^2dx\leq C\int|\na u||\na w_3|^2dx\\
&+\frac{1}{2}\int (H\cdot\na B+B\cdot\na H)\cdot\dot{w}_3dx-\frac{1}{2}\int\na(B\cdot H)\cdot\dot{w}_3dx.
\ea
\ee

Next, we will estimate the three terms on the left of  \eqref{czw6}. For the second term on the left of  \eqref{czw6}, some straightforward calculations  yield  that
\be\la{czw8}\ba
&~~~~\frac{1}{2}\int (H\cdot\na B+B\cdot\na H)\cdot\dot{w}_3dx\\
&=\frac{1}{2}\int (H\cdot\na B+B\cdot\na H)\cdot w_{3t}dx+\frac{1}{2}\int (H\cdot\na B+B\cdot\na H)\cdot(u\cdot\na w_3)dx\\
&=-\frac{1}{2}\int H\cdot\na w_{3t} \cdot Bdx+\frac{1}{2}\int H\cdot\na B\cdot(u\cdot\na w_3)dx\\
&~~~~-\frac{1}{2}\int B\cdot\na w_{3t} \cdot Hdx+\frac{1}{2}\int B\cdot\na H\cdot(u\cdot\na w_3)dx\\
&=-\frac{1}{2}\frac{d}{dt}\left(\int H\cdot\na w_{3} \cdot Bdx+\int B\cdot\na w_{3} \cdot Hdx\right)\\
&~~~~+\frac{1}{2}\int H_t\cdot\na w_{3} \cdot Bdx+\frac{1}{2}\int H\cdot\na w_{3} \cdot B_tdx+\frac{1}{2}\int H\cdot\na B\cdot(u\cdot\na w_3)dx\\
&~~~~+\frac{1}{2}\int B_t\cdot\na w_{3} \cdot Hdx+\frac{1}{2}\int B\cdot\na w_{3} \cdot H_tdx+\frac{1}{2}\int B\cdot\na H\cdot(u\cdot\na w_3)dx\\
&=\frac{d}{dt}R_1+\sum_{k=1}^{6}\hat{I}_k.
\ea
\ee

For $v=(v^1,v^2) ,$ denoting
\be\la{czw7}\ba
\bar{D}(u,v)\triangleq v\cdot \na u-v\div u+\nu \Delta v,~~~
\ea
\ee
we deduce  from   $(\ref{a1})_3$ after integration by parts that
\be\la{czw9}\ba
\hat{I}_1&=-\frac{1}{2}\int (u\cdot\na H)\cdot\na w_{3} \cdot Bdx+\frac{1}{2}\int \bar{D}(u,H)\cdot\na w_{3} \cdot Bdx\\
&=\frac{1}{2}\int \div u  H \cdot\na w_{3} \cdot Bdx+\frac{1}{2}\int  u^iH^j(\p_i\p_jw_{3}^k)B^kdx\\
&~~~~+\frac{1}{2}\int (H\cdot\na w_{3})\cdot(u \cdot \na B)dx+\frac{1}{2}\int \bar{D}(u,H)\cdot\na w_{3} \cdot Bdx,\\[2mm]
\hat{I}_2&=-\frac{1}{2}\int (H\cdot\na w_{3}) \cdot(u\cdot\na B)dx+\frac{1}{2}\int (H\cdot\na w_{3})\cdot\bar{D}(u,B)dx,\\[2mm]
\hat{I}_3&=-\frac{1}{2}\int H^k B^j \p_ku^i \p_iw_3^jdx-\frac{1}{2}\int H^k B^j u^i (\p_k\p_iw_3^j)dx,
\ea
\ee
which implies that
\be\la{czw12}\ba
\sum_{k=1}^3\hat{I}&=\frac{1}{2}\int \div u  H \cdot\na w_{3} \cdot Bdx-\frac{1}{2}\int H \cdot\na u\cdot\na w_{3} \cdot Bdx\\
&~~+\frac{1}{2}\int \bar{D}(u,H)\cdot\na w_{3} \cdot Bdx+\frac{1}{2}\int (H\cdot\na w_{3})\cdot\bar{D}(u,B)dx.
\ea
\ee
Similarly, we also have
\be\la{czw16}\ba
\sum_{k=4}^6\hat{I}&=\frac{1}{2}\int \div u  B \cdot\na w_{3} \cdot Hdx-\frac{1}{2}\int B \cdot\na u\cdot\na w_{3} \cdot Hdx\\
&~~+\frac{1}{2}\int \bar{D}(u,B)\cdot\na w_{3} \cdot Hdx+\frac{1}{2}\int (B\cdot\na w_{3})\cdot\bar{D}(u,H)dx.
\ea
\ee
Following the same arguments as \eqref{czw8}-\eqref{czw12}, it holds that
\be\la{czw17}\ba
&-\frac{1}{2}\int\na(B\cdot H)\cdot\dot{w}_3dx \\
&=\frac{1}{2}\frac{d}{dt}\int(B\cdot H)  \div w_{3}dx-\frac{1}{2}\int(B\cdot H_t)  \div w_{3}dx\\
&\qquad-\frac{1}{2}\int(B_t\cdot H) \div w_{3}dx+\frac{1}{2}\int (B\cdot H)\div(u\cdot\na w_3)dx\\
&=\frac{d}{dt}R_2-\frac{1}{2}\int \div u(H\cdot B) \div w_{3}dx-\frac{1}{2}\int \bar{D}(u\cdot H)\cdot B  \div w_{3}dx\\
&\qquad -\frac{1}{2}\int \bar{D}(u,B)\cdot H  \div w_{3}dx+\frac{1}{2}\int (H\cdot B) (\na u\cdot \na w_{3})dx.
\ea
\ee

Putting  \eqref{czw12}, \eqref{czw16}, \eqref{czw17} into  \eqref{czw6} and using \eqref{czw7}, we obtain that
\be\la{czw21}\ba
&\frac{\mu}{2}\frac{\rm d}{{\rm d}t}\int|\na w_3|^2dx+\int\rho|\dot w_3|^2dx\\
&\le \frac{d}{dt}(R_1+R_2)+C\int|\na u||\na w_3|^2dx+C\int|H||B||\na u||\na w_3|dx\\
&~~+C\int|\na H||B||\na^2 w_3|dx+C\int| H||\na B||\na^2 w_3|dx\\
&\quad +C\int|\na H||\na B||\na  w_3|dx\triangleq \frac{d}{dt}(R_1+R_2)+\sum_{k=1}^5\ti I_k,
\ea
\ee

Noticing that  Gronwall's inequality together with \eqref{bug2} yields that
   \be\la{c11}\sup_{0\le t\le T}\||H|^2\|_{L^2}^2+\int_0^T \||H||\na H| \|_{L^2}^2 dt\leq C \|H_0\|_{L^4}^4\leq C(M),
     \ee we get after direct calculations   that
\be\la{czw22}\ba
R_1+R_2&\le C\int| H|| B||\na  w_3|dx \le \frac{\mu}{4}\|\na w_3\|_{L^2}^2+C\|| H|| B|\|_{L^2}^2,
\ea
\ee
and
\be\la{czw23'}\ba
\tilde{I}_1+\tilde{I}_2 &\le C\int| \na u||\na  w_3|^2dx+ C\int| H|^2| B|^2|\na u|dx\\
&\le C\|\na u\|_{L^2} \|\na w_3\|_{L^4}^2+C\||H||B|\|_{L^4}^2\|\na u\|_{L^2} \\
&\le \varepsilon\|\na^2 w_3\|_{L^2}^2+C(\ve)\|\na u\|_{L^2}^2\|\na w_3\|_{L^2}^2+C\|\na u\|_{L^2}^2\||H||B|\|_{L^2}^2\\
&\quad+C\||B||\na H|\|_{L^2}^2+C\||H||\na B|\|_{L^2}^2,\ea\ee \be\ba
\tilde{I}_3+\tilde{I}_4&\le \varepsilon\|\na^2 w_3\|_{L^2}^2+C(\ve)\||B||\na H|\|_{L^2}^2+C(\ve)\||H||\na B|\|_{L^2}^2,\ea\ee \be \la{czw23''}\ba
\tilde{I}_5&\le C\|\na H\|_{L^2} \|\na B\|_{L^4}\|\na w_3\|_{L^4}\\
&\le C\|\na H\|_{L^2} \|\na w_3\|_{L^4}^2+C\|\na H\|_{L^2} \|\na B\|_{L^4}^2 \\
&\le \varepsilon\|\na^2 w_3\|_{L^2}^2+C(\ve)\|\na H\|_{L^2}^2\|\na w_3\|_{L^2}^2\\
&\quad+C\|\na H\|_{L^2}^2\|\na B\|_{L^2}^2+ C\|\na^2 B\|_{L^4}^2.
\ea
\ee
Furthermore, applying the standard $L^p$-estimate to \eqref{sasw3} gives
\be\la{lbq-jia43}\ba
 \|\na^2 w_3\|_{L^2}^2 &\leq  C\|\n^{1/2}\dot{\omega}_3\|_{L^2}^2+C\||H| |\na B|\|_{L^2}^2+C\||B| |\na H|\|_{L^2}^2,
 \ea\ee
which together with \eqref{czw23'}-\eqref{czw23''}  and choosing $\varepsilon$ suitably small yields
\be\la{lvczw23}\ba
\sum_{i=1}^{5}\tilde{I}_i
&\le \frac{1}{2}\|\n^{1/2}\dot{\omega}_3\|_{L^2}^2+C(\|\na u\|_{L^2}^2+\|\na H\|_{L^2}^2)\|\na w_3\|_{L^2}^2\\
&\quad+C\||B||\na H|\|_{L^2}^2+C\||H||\na B|\|_{L^2}^2\\
&\quad+C\|\na u\|_{L^2}^2\||H||B|\|_{L^2}^2+C\|\na H\|_{L^2}^2\|\na B\|_{L^2}^2+ C\|\na^2 B\|_{L^4}^2.
\ea
\ee
Putting \eqref{lvczw23}  into \eqref{czw21} shows
\be\la{czw27}\ba
&\frac{\mu}{2}\frac{\rm d}{{\rm d}t}\int|\na w_3|^2dx+\frac{1}{2}\int\rho|\dot w_3|^2dx\\
&\le\frac{d}{dt}(R_1+R_2)+C(\|\na u\|_{L^2}^2+\|\na H\|_{L^2}^2) \|\na w_3\|_{L^2}^2+C\|\na^2 B\|_{L^2}^2\\
 &\quad+C\|\na u\|_{L^2}^2(\||H||B|\|_{L^2}^2+\|\na B\|_{L^2}^2)+C\||B||\na H|\|_{L^2}^2+C\||H||\na B|\|_{L^2}^2.
\ea
\ee

Furthermore, multiplying \eqref{a1}$_3$ by $H|B|^2$, and multiplying \eqref{a1'} by $B|H|^2$, then adding the upper two resultant equations together and  integrating by parts over $\mathbb{R}^2$,   we have
\be\la{hb5}\ba
&\frac{1}{2}\frac{d}{dt}\||B||H|\|_{L^2}^2+\nu\||\na H||B|\|_{L^2}^2+\nu\||\na B||H|\|_{L^2}^2\\
&\le  C\int |\na H||H|  |\na B| |B| dx+ C\int |\na u||H|^2|B|^2dx\\
&\le \||\na H||B|\|_{L^2}\|\na B\|_{L^4}\|H\|_{L^4}+C \|\na u\|_{L^2} \||H||B|\|_{L^2}  \|\na(|H||B|)\|_{L^2} \\
&\le  \frac{\nu}{2}(\||\na H||B|\|_{L^2}^2+\||\na B||H|\|_{L^2}^2)+C(M) \|\na B\|_{L^4}^2+C \|\na u\|_{L^2}^2 \||H||B|\|_{L^2}^2.
\ea\ee
Similar to \eqref{czh4} and \eqref{czh5}, for $\theta=  0,1 $,  it holds
 \be\la{czh4'}\ba  \sup_{0\le t\le\si( T)} t^\theta\|B\|_{H^1}^2  +\int_0^{\si(T)} t^\theta\|\na
B\|_{H^1}^2 dt\le C\|B_0\|_{H^{1-\theta}}^2. \ea\ee
This combined with Gronwall's inequality, \eqref{hb5} and  \eqref{a16} gives
\be\la{hb6}\ba
& \sup_{0\le t\le \si(T)}(\||B||H|\|_{L^2}^2)+\int_0^{\si(T)} \left(\||\na H||B|\|_{L^2}^2+ \||\na B||H|\|_{L^2}^2\right)dt \\
&\le  C\||H_0||B_0|\|_{L^2}^2+C (M) \int_0^{\si(T)}\left(\|\na B\|_{L^2}^2+\|\Delta B\|_{L^2}^2\right)dt\\
&\leq C\|H_0\|_{L^4}^2(\|B_0\|_{L^2}^2+\|\na B_0\|_{L^2}^2)+ C(M)\|B_0\|_{H^1}^2\leq C(M)\|B_0\|_{H^1}^2
\ea\ee
and
\be\la{hb7}\ba
& \sup_{0\le t\le \si(T)}(t\||B||H|\|_{L^2}^2)+\int_0^{\si(T)} t\left(\||\na H||B|\|_{L^2}^2+ \||\na B||H|\|_{L^2}^2\right)dt \\
&\le  C\int_0^{\si(T)}\||B||H|\|_{L^2}^2dt +C (M) \int_0^{\si(T)}t\left(\|\na B\|_{L^2}^2+\|\Delta B\|_{L^2}^2\right)dt\\
&\le C\int_0^{\si(T)}\|H\|_{L^4}^2 (\|B\|_{L^2}^2+\|\na B \|_{L^2}^2)dt+ C(M)\|B_0\|_{L^2}^2\le C(M)\|B_0\|_{L^2}^2.
\ea\ee

Now, applying Gronwall's inequality to  \eqref{czw27}, together with  \eqref{c11},  \eqref{czw22}, \eqref{czh4'},  \eqref{hb6} and  \eqref{hb7}  shows
\be\la{czw28}\ba
\sup_{0\le t\le \si(T)}\|\na w_3\|_{L^2}^2+\int_0^{\si(T)}\int\rho|\dot w_3|^2dxdt
\le C(M)\|B_0\|_{H^1}^2,
\ea
\ee
   \be\la{czw29} \ba \sup_{0\le t\le \si(T)}t\|\na w_3\|_{L^2}^2+\int_0^{\si(T)}t\int\rho|\dot w_3|^2dxdt
\le C(M)\|B_0\|_{H^0}^2.\ea\ee
Since the solution operator $B_0\mapsto B$ and $B\mapsto w_3$ is linear, so that $B_0\mapsto w_3$ is also linear, we thus conclude from
\eqref{czw28}  and \eqref{czw29} in a manner similar to the derivation of  \eqref{lbqnew-czh} that
   \be\la{czw30} \ba \sup_{0\le t\le \si(T)}t^{1-\theta}\|\na w_3\|_{L^2}^2+\int_0^{\si(T)}t^{1-\theta}\int\rho|\dot w_3|^2dxdt
\le C(M)\|H_0\|_{\dot{H}^\theta}^2 . \ea\ee

Finally, choosing $w_{10}=u_0$ and $B_0= H_0$, so that $w_1+w_2+w_3=u$ and $B=H$, we immediately obtain \eqref{uv1} from \eqref{uu4}, \eqref{uu3}  and \eqref{czw30}.  Thus, we finish the proof of Lemma \ref{zc1}.

We now proceed to derive a uniform (in time) upper bound for the
density, which turns out to be the key to obtain all the higher
order estimates.
We will use an approach motivated by   \cite{lx,hlx1}).

\begin{lemma}\la{le7}
There exists a positive constant
   $\ve_0=\ve_0 (\on ,M) $
    depending    on  $\mu ,  \lambda ,  \nu,  \ga ,  a ,  \on, $ $\beta, $ $ N_0,$ and $M$  such that,
    if  $(\n,u,H)$ is a smooth solution  of
   (\ref{a1})-(\ref{h2})     on $\O \times (0,T] $
   satisfying (\ref{z1}) and the assumptions in Theorem \ref{th1}, then
      \be\la{lv102}\sup_{0\le t\le T}\|\n(t)\|_{L^\infty}  \le
\frac{7\bar \n }{4}  ,\ee
      provided $C_0\le \ve_0 . $

   \end{lemma}

{\it Proof.}
  First, we rewrite the equation of the mass conservation
$(\ref{a1})_1$ as \be \la{z.3} D_t \n=g(\n)+b'(t), \ee where \bnn
D_t\n\triangleq\n_t+u \cdot\nabla \n ,\quad
g(\n)\triangleq-\frac{\n^{\ga+1 }}{2\mu+\lambda}  ,
\quad b(t)\triangleq-\frac{1}{2\mu+\lambda} \int_0^t\n \left(F+\frac{1}{2} |H|^2\right)dt. \enn

Next, it follows from \eqref{h19},   \eqref{uq2}, and \eqref{z.1}   that for $t>0$ and $p\in [2,\infty),$
 \be\la{z.2}\ba   \|\na \left(F+\frac{1}{2}|H|^2\right)\|_{L^p}&\le C\|\na F\|_{L^p}+C\|\na |H|^2\|_{L^p}\\
  &\le C(p)\|\n \dot u \|_{L^p} +C(p)\||H||\na H|\|_{L^p}\\
 &\le C(p,\on,M)(1+t)^5\left(\|\n^{1/2}\dot u \|_{L^2}+\|\na\dot u \|_{L^2}\right)\\
 &\quad+C(p)\|\na H\|_{L^2}\left(\|H\|_{L^2}+\|\triangle H\|_{L^2}\right)
 \ea\ee
 where in the last inequality we have used  following simple facts:
  \be\la{AMSS3.1}\ba \||H||\na H|\|_{L^p}
 &\le \|H\|_{L^{2p}} \|\na H\|_{L^{2p}}\le C(p)\|H\|_{L^2}^{\frac{1}{p}}\|\na H\|_{L^2}\|\triangle H\|_{L^2}^{\frac{p-1}{p}}\\
 &\le C(p)\|\na H\|_{L^2}\left(\|H\|_{L^2}+\|\triangle H\|_{L^2}\right) \ea\ee
 and
 \be \la{o3.7} \sup_{0\le s\le t}\int\n  (1+|x|^2)^{1/2}  dx\le  C(M)(1+t) \ee
 which is obtained after multiplying $\eqref{a1}_1$ by $(1+|x|^2)^{1/2} $ and integrating the resulting equality over $\O$ by parts(see \cite[(3.40)]{lx1} also).
Choosing $q=2$ in the Gagliardo-Nirenberg inequality  \eqref{g2}  and using \eqref{z.2}, we deduce that
for $r\triangleq 4+4/\beta$ and $\de_0\triangleq (2r+(1-\beta)(r-2))/(3r-4)\in(0,1),$
\be\la{lbqnew-jia32}\ba
 |b(\si(T)) | 
 &\le C(\on)\int_0^{\si(T)} \si^{-\frac{2r+(1-\beta)(r-2)}{4(r-1)}}\left(\si^{1-\beta}(\|F\|_{L^2}^2+\||H|^2\|_{L^2}^2)\right)
^{\frac{r-2}{4(r-1)}}\\
&\quad \cdot \left(\si^2\|\na  (F+|H|^2) \|^{2}_{L^r}\right)^\frac{r}{4(r-1)}dt\\
&\le  C(\on,M) \int_0^{\si(T)} \si^{-\frac{2r+(1-\beta)(r-2)}{4(r-1)}}\left(\si^2\|\na  (F+|H|^2)\|^{2}_{L^r} \right)^\frac{r}{4(r-1)}dt\\
&\le  C(\on,M) \left(\int_0^{\si(T)}
\si^{-\de_0}dt\right)^\frac{3r-4}{4(r-1)}  \left(\int_0^{\si^(T)} \si^2\|\na  (F+|H|^2)\|^{2}_{L^r}  dt\right)^\frac{r}{4(r-1)}\\
&\le  C(\on,M) C_0^\frac{r}{4(r-1)} \ea\ee
where in the second inequality we have used  \eqref{lbqnew-czh}  and \eqref{uv1},
and in the last   one  we used the following estimate  which comes from \eqref{a16} and \eqref{h27}
\be\la{lbqnew-jia33}\ba
 &\int_0^{\si^(T)} \si^2\|\na  (F+|H|^2)\|^{2}_{L^r}  dt\\
 &\le   C\int_0^{\si^(T)} \si^2\left(\|\n^{1/2}\dot u \|_{L^2}^2+\|\na\dot u \|_{L^2}^2\right) dt\\&\quad + C\int_0^{\si^(T)} \si^2\|\na H\|_{L^2}^2(\|H\|_{L^2}^2+\|\triangle H\|_{L^2}^2)dt\\
 &\le CC_0+C \left(\sup_{t\in
[0,\si(T)]}(\si \|\na H\|_{L^2}^2)\right)\int_0^{\si^(T)}(\|H\|_{L^2}^2+\si \|\triangle H\|_{L^2}^2)dt\\
&\le  CC_0.
 \ea\ee
  Hence, \eqref{lbqnew-jia32}  combined with \eqref{z.3} yields that
   \be\la{a103}\sup_{t\in
[0,\si(T)]}\|\n\|_{L^\infty} \le \on
+C(\on,M)C_0^{1/4}\le\frac{3 \bar\n  }{2},\ee
 provided $$C_0\le \ve_1\triangleq\min\{1, (\on/(2C(\on,M)))^{4}\}. $$

Furthermore, it follows from  \eqref{h19} and  \eqref{ly8}   that for $t\in[\si(T),T],$
\be \la{hg2} \ba \|F+|H|^2\|_{H^1} 
& \le C \left(\|\na u \|_{L^2}+\|P \|_{L^2}+\||H|^2\|_{L^2}\right)\\
&~~~~+C \left(\|\n  \dot u \|_{L^2}+\||H||\na H|\|_{L^2}+\|\na |H|^2\|_{L^2}\right)\\
& \le C( \on) C_0^{1/2}  t^{-1/2}, \ea\ee
which together with \eqref{g2} and  \eqref{z.2} shows
\be\la{hg1}\ba  \int_{\si(T)}^T\|F+|H|^2\|_{L^\infty}^4dt
  &\le C\int_{\si(T)}^T\|F+|H|^2\|_{L^{72}}^{\frac{35}{9}}\|\na F+\na |H|^2\|_{L^{72}}^{\frac{1}{9}}dt\\
&\le C(\on,M)C_0^{35/18}\int_{\si(T)}^T \left[ t^{-\frac{25}{18}}(\|\n^{1/2}\dot u\|_{L^2} +\|\na\dot u\|_{L^2} )^{\frac{1}{9}}\right.\\
&\qquad~+\left. t^{-\frac{35}{18}}(\|\na H\|_{L^2}(\|H\|_{L^2}+\|\triangle H\|_{L^2}))^{\frac{1}{9}}\right]dt\\
&\le C(\on,M)C_0^{35/18},\ea\ee
where in the last inequality, one has used \eqref{z1} and \eqref{a16}.
This shows that for all $\si(T)\le t_1\le t_2\le T,$ \bnn\ba  |b(t_2)-b(t_1)|    &\le  C(\on) \int_{t_1}^{t_2} \|F+|H|^2\|_{L^\infty}dt  \\&\le \frac{1}{2\mu+\lambda}(t_2-t_1)+ C(\on,M) \int_{\si(T)}^T \|F+|H|^2\|_{L^\infty}^4dt\\ &\le \frac{1}{2\mu+\lambda}(t_2-t_1)+  C(\on,M) C_0^{35/18}   ,  \ea\enn which implies that  one can
choose $N_1$ and $N_0$ in (\ref{a100}) as: \bnn
N_1=\frac{ 1}{2\mu+\lambda},\quad  N_0=C(\on,M)C_0^{35/18} .\enn  Hence, we set $\bar\zeta= 1 $ in (\ref{a101}) since for all $  \zeta \ge
 1,$
$$ g(\zeta)=-\frac{ \zeta^{\ga+1}}{2\mu+\lambda} \le -N_1
=- \frac{ 1}{2\mu+\lambda}.  $$   Lemma
\ref{le1} and (\ref{a103}) thus lead to \be\la{a102} \sup_{t\in
[\si(T),T]}\|\n\|_{L^\infty}\le  \frac{ 3\bar \n
}{2}  +N_0 \le
\frac{7\bar \n }{4} ,\ee provided $$ C_0\le
\ve_0\triangleq\min\{\ve_1,\ve_2 \}, \quad\mbox{ for
}\ve_2\triangleq \left(\frac{ \bar \n }{4C(\on,M) }\right)^{18/35}.$$
The combination of (\ref{a103}) with (\ref{a102}) completes the
proof of Lemma \ref{le7}.

With Lemma \ref{le5} and Lemma \ref{le7} at hand, we are now in a position to prove Proposition \ref{pr1}.

{\it Proof of Proposition  \ref{pr1}.} It follows from \eqref{h27} that
  \be\la{end-h27}
  A_1(T)+A_2(T)+\int_0^T\si\|P\|_{L^2}^2dt\le  C_0^{1/2}
  \ee
provided
$$C_0\le\varepsilon_3\triangleq(C(\bar{\n}))^{-2}.$$

It remains to estimate $A_3(\si(T))$. Indeed, using \eqref{end-h27}, \eqref{lbqnew-czh}, and \eqref{uv1},  it holds that
  \be\la{end-A3}\ba
A_3(\si(T)) \leq& \sup_{  0\le t\le \si(T)   }\left(\left(\sigma^{1-\beta} \|\na u\|_{L^2}^2\right)^{(2\beta+1)/(4\beta)}  \left(\sigma  \|\na u\|_{L^2}^2\right)^{(2\beta-1)/(4\beta)}\right)\\
& +\sup_{  0\le t\le \si(T)   }\left(\left(\sigma^{1-\beta} \|\na H\|_{L^2}^2\right)^{(2\beta+1)/(4\beta)}
 \left(\sigma  \|\na H\|_{L^2}^2\right)^{(2\beta-1)/(4\beta)}\right)\\
\leq& C(\bar{\n},M)A_1^{(2\beta-1)/(4\beta)}(T)\leq C(\bar{\n},M)C_0^{(2\beta-1)/(8\beta)}\leq  C_0^{\delta_0}
\ea
\ee
provided
$$C_0\le\varepsilon_4\triangleq   C(\bar{\n},M)^{(-72\beta)/(2\beta-1)}.$$

Letting $\varepsilon\triangleq\min\{\varepsilon_0, \varepsilon_3, \varepsilon_4\}$, we obtain \eqref{z2} directly from \eqref{lv102}, \eqref{end-h27} and \eqref{end-A3} and finish the proof of Proposition  \ref{pr1}.

\section{\la{se5} A priori estimates (II): higher order estimates }

From now on, for smooth initial data $(\n_0,u_0,H_0)$  satisfying \eqref{co1} and \eqref{h7},   assume that
  $(\n,u,H)$ is a smooth solution of (\ref{a1})-(\ref{h2}) on $\O \times (0,T] $ satisfying (\ref{z1}).
Then, we  derive some necessary uniform  estimates on the spatial gradient of
the smooth solution $(\n,u,H)$.

\begin{lemma}\la{le4}     There is a positive constant $C $ depending only on $T,\mu ,  \lambda , \nu,  \ga ,  a ,  \on,$ $ \beta, N_0, $ $M,q,$   and $\|\n_0\|_{H^1\cap W^{1,q}} $ such that
  \be\la{pa1}\ba
  &\sup_{0\le t\le T}\left(\norm[H^1\cap W^{1,q}]{ \rho} +\|\nabla   u\|_{L^2}+\|\nabla  H\|_{L^2} +  t\|\na^2 u\|^2_{L^2}  \right)\\&+\int_0^T \left(\|\nabla^2  u\|_{L^2}^2+\|\na^2 H\|_{L^2}^2+\|\nabla^2 u\|_{L^q}^{(q+1)/q}+t\|\nabla^2 u\|_{L^q}^2 \right)dt\le C .\ea
  \ee
\end{lemma}

{\it Proof.}  First, it follows from \eqref{h29}, \eqref{ly12},  Gronwall's inequality, and \eqref{a16} that
     \be\la{a93}\ba&
  \sup_{t\in[0,T]}\left(\|\nabla u\|_{L^2}^2+\|\nabla H\|_{L^2}^2\right)\\& + \int_0^{T} \left(\|\n^{1/2}\dot{u}\|_{L^2}^2+\|\Delta H\|_{L^2}^2+\||H||\nabla H|\|_{L^2}^2\right)dt
  \le C,\ea
  \ee which together with \eqref{h18}   shows \be\la{hc1} \int_0^T\left(\|\na u\|_{L^4}^4+\|\na H\|_{L^4}^4\right)dt\le C.\ee
This combined with \eqref{nlv4}  shows
\be\la{b19-1}\ba
\sup_{0\le t\le T}(t\||H||\na H|\|_{L^2}^2)+\int_0^Tt(\|\triangle |H|^2\|_{L^2}^2+\||H||\triangle H|\|_{L^2}^2)dt\le C.
\ea\ee
Multiplying \eqref{lv3.46} by $t$ and integrating the resulting inequality over $(0,T)$ combined with  \eqref{a93}, \eqref{hc1} and \eqref{b19-1}  lead to
\be\la{b19}\ba
& \sup_{0\le t\le T}t ( \|\n^{1/2}\dot{u}\|_{L^2}^2 + \||H||\na H|\|_{L^2}^2)\\
  &~~~~+\int_0^Tt(\|\na\dot{u}\|_{L^2}^2+\|\triangle |H|^2\|_{L^2}^2+\||H||\triangle H|\|_{L^2}^2)dt\le C.
\ea\ee

Next, we prove (\ref{pa1}) by using Lemma \ref{le9} as in
\cite{hlx}. For $   p\in [2, q],$ $|\nabla\n|^p$ satisfies \bnnn \ba
& (|\nabla\n|^p)_t + \text{div}(|\nabla\n|^pu)+ (p-1)|\nabla\n|^p\text{div}u  \\
 &+ p|\nabla\n|^{p-2}(\nabla\n)^t \nabla u (\nabla\n) +
p\n|\nabla\n|^{p-2}\nabla\n\cdot\nabla\text{div}u = 0.\ea
\ennn
Thus, \be\la{L11}\ba
\frac{d}{dt} \norm[L^p]{\nabla\n}  &\le
 C(1+\norm[L^{\infty}]{\nabla u} )
\norm[L^p]{\nabla\n} +C\|\na^2u\|_{L^p}\\ &\le
 C(1+\norm[L^{\infty}]{\nabla u} )
\norm[L^p]{\nabla\n} +C\|\n\dot u\|_{L^p}+C\||H||\na H|\|_{L^p}, \ea\ee
due to
\be
\la{ua1}\|\na^2 u\|_{L^p}\le   C\left(\|\n\dot u\|_{L^p}+ \|\nabla
P \|_{L^p}+\||H||\na H|\|_{L^p}\right),\ee
which follows from the standard
$L^p$-estimate for the following elliptic system:
 \bnn  \mu\Delta
u+(\mu+\lambda)\na {\rm div}u=\n \dot u+\na P+\frac{1}{2}\na|H|^2-\text{div}(H\otimes H),\quad \, u\rightarrow
0\,\,\mbox{ as } |x|\rightarrow \infty. \enn

Next, it follows from
  the  Gargliardo-Nirenberg inequality, \eqref{a93},  and  \eqref{h19}   that
 \be\la{419}\ba  &\|\div u\|_{L^\infty}+\|\o\|_{L^\infty}  \\
 &\le C \|F\|_{L^\infty}+C\|P\|_{L^\infty}+C\||H|^2\|_{L^\infty} +C\|\o\|_{L^\infty}\\ &\le C(q) +C(q) \|\na F\|_{L^q}^{q/(2(q-1))}+C(q) \|\na |H |^2\|_{L^q}^{q/(2(q-1))} +C(q) \|\na \o\|_{L^q}^{q/(2(q-1))}\\ &\le C(q) +C(q) \left(\|\n\dot u\|_{L^q}+\||H||\na H|\|_{L^q}\right)^{q/(2(q-1))} , \ea\ee
 which, together with
Lemma \ref{le9}, \eqref{ua1} and \eqref{a93}, yields that
    \be\la{b24}\ba   \|\na
u\|_{L^\infty }\le& C\left(\|{\rm div}u\|_{L^\infty }+
\|\o\|_{L^\infty } \right)\log(e+\|\na^2 u\|_{L^q}) +C\|\na
u\|_{L^2} +C \\\le& C\left(1+\|\n\dot u\|_{L^q}^{q/(2(q-1))}+\||H||\na H|\|_{L^q}^{q/(2(q-1))}\right)\\
&\cdot\log(e+\|\rho \dot u\|_{L^q} +\||H||\na H| \|_{L^q}+\|\na \rho\|_{L^q})\\
\le& C\left(1+\|\n\dot u\|_{L^q} +\||H||\na H|\|_{L^q}\right)\log(e+   \|\na \rho\|_{L^q}) . \ea\ee

Next, it follows from the H\"older inequality   and \eqref{z.2}   that
 \bnn\la{b22}\ba
 \| \rho \dot u\|_{L^q} & \le
  \| \rho \dot u\|_{L^2}^{2(q-1)/(q^2-2)}\|\n\dot u\|_{L^{q^2}}^{q(q-2)/(q^2-2)}\\ & \le
 C \| \rho \dot u\|_{L^2}^{2(q-1)/(q^2-2)}\left(\| \rho^{1/2} \dot u\|_{L^2}+\|\na\dot u\|_{L^2}\right)^{q(q-2)/(q^2-2)}\\ & \le
 C \| \rho^{1/2}  \dot u\|_{L^2} +C \| \rho^{1/2} \dot u\|_{L^2}^{2(q-1)/(q^2-2)}\|\na \dot u\|_{L^2}^{q(q-2)/(q^2-2)}   , \ea\enn which combined with   \eqref{a93} and  \eqref{b19} implies that\be\la{4a2}   \ba &\int_0^T\left(\| \rho \dot u\|_{L^q}^{1+1 /q}+t\| \rho \dot u\|_{L^q}^2\right)  dt\\ & \le C   \int_0^T\left( \| \rho^{1/2}  \dot u\|_{L^2}^2+ t\|\na \dot u\|_{L^2}^2 + t^{-(q^3-q^2-2q-1)/(q^3-q^2-2q )}\right)dt+C \\
 &\le  C.\ea\ee
 Moreover, we have by \eqref{AMSS3.1},  \eqref{a16}  and \eqref{a93}  that
 \be\la{lbq-jia8}\ba
&\int_0^T \left(\||H||\na H|\|_{L^q}^{1+1/q}+ \||H||\na H|\|_{L^q}^2\right)dt\\
&\le C\int_0^T \left(\|\na^2 H\|_{L^{2}}^{1-1/q^2}+ \|\na^2 H\|_{L^{2}}^{2-2/q}\right)dt\\
&\le C\int_0^T \left(1+ \|\na^2 H\|_{L^{2}}^2\right)dt\le C
\ea\ee

Then,
 substituting \eqref{b24} into \eqref{L11} where $p=q$, we deduce from Gronwall's inequality, \eqref{a93} and \eqref{4a2}-\eqref{lbq-jia8}  that \be \la{b30} \sup\limits_{0\le t\le T}\|\nabla
\rho\|_{L^q}\le C,\ee  which, along with \eqref{ua1}, \eqref{4a2} and \eqref{lbq-jia8}, shows

\be\la{mhd7}
\int_0^T\left(\|\na^2u\|_{L^q}^{(q+1)/q}+t\|\na^2u\|_{L^q}^2\right)dt\leq C.
\ee

Finally,
taking $p=2$ in (\ref{L11}), one gets by using (\ref{a93}), \eqref{mhd7}, 
 and Gronwall's inequality that \bnn \sup\limits_{0\le
t\le T}\|\nabla \n\|_{L^2}\le C,\enn which, together with (\ref{b30}),  (\ref{a93}),
(\ref{ua1}),   (\ref{b19}),and
(\ref{mhd7}), yields (\ref{pa1}). The proof of Lemma \ref{le4} is
completed.

Next, we will show the following spatial weighted mean estimate of the density, which has been proved in \cite[Lemma 4.2]{lx1}.
\begin{lemma} \la{le6} There is a positive  constant $C $ depending only on  $T,\mu ,  \lambda ,  \nu,  \ga ,  a ,  \on, \beta, N_0,$ $  M,q,$    and $\|\na(\bar x^a\n_0)\|_{L^2\cap L^q}  $ such that
  \be \la{q} \ba
  &\sup_{0\le t\le T}   \|  \bar x^a\n \|_{L^1\cap H^1\cap W^{1,q}}  \le C .\ea
  \ee
\end{lemma}

\begin{lemma}\la{newle}
There exists a positive  constant $C $ depending only on  $T,\mu ,  \lambda , \nu,  \ga ,  a ,  \on, \beta, N_0,$ $ M,q,$    and $\||H_0|^2\bar{x}^a\|_{L^1}$ such that
\be
\ba\la{gj10}
\sup_{0\le t\le T} \|H\bar{x}^{a/2}\|_{L^2}^2 +\int_{0}^{T} \|\na H\bar{x}^{a/2}\|_{L^2}^2dt\le C,
\ea
\ee
\be
\ba\la{gj10'}
\sup_{0\le t\le T} \left(t\|\na H\bar{x}^{a/2}\|_{L^2}^2\right)+\int_{0}^{T} t\|\Delta H\bar{x}^{a/2}\|_{L^2}^2dt\le C.
\ea
\ee
\end{lemma}

{\it Proof.}
First, it follows from \eqref{ljo1} and  \eqref{a93} that for any $\eta\in(0,1]$ and any $s>2,$
 \be\la{a4.21} \|u\bar x^{-\eta}\|_{L^{s/\eta}}\le C(\eta,s).\ee
Similar to the proof of \eqref{lbqnew-gj10} (~or \cite[Lemma 4.1]{chl}), multiplying  \eqref{a1}$_3$ by $H\bar{x}^a$, integration by parts together with \eqref{a93} and   \eqref{a4.21} yields
\be\la{mhd8}\ba
&\frac{1}{2}\left(\|H\bar{x}^{a/2}\|_{L^2}^2\right)_t+\nu \|\na H\bar{x}^{a/2}\|_{L^2}^2\le  C\|H\bar{x}^{a/2}\|_{L^2}^2+ \frac{\nu}{2}\| \na H \bar{x}^{a/2}\|_{L^2}^2,
\ea\ee
which together with Gronwall's inequality yields \eqref{gj10}.

Now, multiplying  \eqref{a1}$_3$ by $\Delta H\bar{x}^a$, integrating the resultant equation by parts over $\mr^2$, it follows from the similar arguments as \eqref{mm1} that
\be\la{AMSS5}\ba
&\frac{1}{2}\left(\int |\na H|^2\bar{x}^adx\right)_t+\nu \int |\Delta H|^2\bar{x}^adx\\
\le& C\int|\na H| |H| |\na u| |\na\bar{x}^a|dx+C\int|\na H|^2|u| |\na\bar{x}^a|dx+C\int|\na H| |\Delta H| \bar{x}^adx\\
&+C\int |H||\na u||\Delta H|\bar{x}^adx+C\int |\na u||\na H|^2  \bar{x}^adx
\triangleq \sum_{i=1}^5 J_i.
\ea\ee
Using Gagliardo-Nirenberg inequality, \eqref{pa1}, \eqref{gj10} and \eqref{a4.21}, it holds
\be\la{AMSS6}\ba
J_1\le &C\int|\na H| |H| |\na u| \bar{x}^{a} (\bar{x}^{-1} |\na\bar{x}|)dx\\
\le & C\|H\bar{x}^{a/2}\|_{L^4}^4+C\|\na u\|_{L^4}^4+C\|\na H\bar{x}^{a/2}\|_{L^2}^2\\
\le & C\|H\bar{x}^{a/2}\|_{L^2}^2\left(\|\na H\bar{x}^{a/2}\|_{L^2}^2+\|H\bar{x}^{a/2}\|_{L^2}^2\right)+C\|\na u\|_{L^4}^4+C\|\na H\bar{x}^{a/2}\|_{L^2}^2\\
\le &C+C\|\na^2 u\|_{L^2}^2+C\|\na H\bar{x}^{a/2}\|_{L^2}^2,
\ea\ee
\be\la{AMSS7}\ba
J_2
&\leq C\int |\na H|^{(2a-1)/a}\bar{x}^{(2a-1)/2}|\na H|^{1/a}|u|\bar{x}^{-1/4} \bar{x}^{-1/4} |\na\bar{x}|dx\\
&\leq  C\||\na H|^{(2a-1)/a}\bar{x}^{(2a-1)/2}\|_{L^{\frac{2a}{2a-1}}} \|u\bar{x}^{-1/4}\|_{L^{4a}}\||\na H|^{1/a}\|_{L^{4a}}\\
&\leq  C\|\na H \bar{x}^{a/2} \|_{L^2}^2+  C\|\na H \|_{L^4}^2\\
&\leq  C\|\na H \bar{x}^{a/2}\|_{L^2}^2+  \varepsilon\|\Delta H \bar{x}^{a/2}\|_{L^2}^2,
\ea\ee
\be\la{AMSS8}\ba
J_3+J_4\le &\varepsilon\|\Delta H\bar{x}^{a/2}\|_{L^2}^2+C\|\na H\bar{x}^{a/2}\|_{L^2}^2+C\|H\bar{x}^{a/2}\|_{L^4}^4+C\|\na u\|_{L^4}^4\\
\le &\varepsilon\|\Delta H\bar{x}^{a/2}\|_{L^2}^2+C+C\|\na H\bar{x}^{a/2}\|_{L^2}^2+C\|\na^2 u\|_{L^2}^2,
\ea\ee
\be\la{AMSS9}\ba
J_5\le &C\|\na u\|_{L^\infty} \|\na H\bar{x}^{a/2}\|_{L^2}^2\\
\le &C \|\na u\|_{L^2}^{(q-2)/(2q-2)}\|\na^2 u\|_{L^q}^{q/(2q-2)} \|\na H\bar{x}^{a/2}\|_{L^2}^2\\
\le &C(1+\|\na^2 u\|_{L^q}^{(q+1)/q})\|\na H\bar{x}^{a/2}\|_{L^2}^2.
\ea\ee
Submitting \eqref{AMSS6}-\eqref{AMSS9} into \eqref{AMSS5} and choosing $\varepsilon$ suitably small, we have
\be\la{AMSS10}\ba
 &\frac{1}{2}\left(\int |\na H|^2\bar{x}^adx\right)_t+\nu \int |\Delta H|^2\bar{x}^adx\\
&\quad\le  C(1+\|\na^2 u\|_{L^q}^{(q+1)/q})\|\na H\bar{x}^{a/2}\|_{L^2}^2+C(\|\na^2 u\|_{L^2}^2+1),
\ea\ee
which multiplied by $t$, then together with  Gronwall's inequality,  \eqref{gj10}  and \eqref{pa1} yields  \eqref{gj10'}.
The proof of Lemma \ref{newle} is finished.

\begin{lemma}\label{lem4.5v}  There is a positive  constant $C $ depending only on  $T,\mu ,  \lambda , \nu,  \ga ,  a ,  \on, \beta, N_0,M,q,$  $\||H_0|^2\bar{x}^a\|_{L^1}$ and $\|\na(\bar x^a\n_0)\|_{L^2\cap L^q}  $ such that
\begin{equation}\la{gj13}\ba
 \sup_{0\leq t\leq T }t\left(\|\n^{1/2}u_t\|^2_{L^2}+\|H_t\|^2_{L^2} \right)+\int_0^Tt\|\na u_t\|_{L^2}^2+t\|\na H_t\|_{L^2}^2dt\le C,\ea
\end{equation}
 \be\la{gj13'} \sup_{0\leq t\leq T }t\left(\|\na u\|^2_{H^1}+\|\na H\|^2_{H^1} \right)\le C. \ee
\end{lemma}

{\it Proof.}
First, the combination of \eqref{a4.21} with \eqref{q} gives  that
for any $\eta\in(0,1]$ and any $s>2,$
 \be\la{5.d2}\|\n^\eta u \|_{L^{s/\eta}}+ \|u\bar x^{-\eta}\|_{L^{s/\eta}}\le C(\eta,s).\ee

Multiplying    equations   \eqref{a1}$_2$  by  $ u_t$ and    integrating by parts, with the similar arguments as  the proof of  \eqref{m1} (or \cite[Lemma 3.2]{lx}), we have by \eqref{a95}, \eqref{5.d2}, \eqref{pa1}  and  \eqref{q} that
\be\la{3r1}\ba& \frac{d}{dt}\int\left(( \mu+\lm)(\div u)^2+\mu|\na u|^2\right)dx+\int \n|  u_t|^2dx\\
&\le 2\int P\div u_tdx+C\int\n|u|^2|\na u|^2dx+ \int\left(H\cdot \na H-\frac{1}{2}\na |H|^2\right)u_tdx\\
&\le 2\frac{d}{dt}\int P\div udx+ C\|\na^2u\|_{L^2}^2+C-\int H\cdot \na u_t\cdot H-\frac{1}{2}  |H|^2 \div u_tdx\\
&\le \frac{d}{dt}\Psi(t)+ C\|\na^2u\|_{L^2}^2+C+\int |H_t||H||\na u|dx\\
&\le \frac{d}{dt}\Psi(t)+ C\|\na^2u\|_{L^2}^2+C+\varepsilon\|H_t\|_{L^2}^2+C\|H\|_{L^4}^2\|\na u\|_{L^4}^2\\
&\le \frac{d}{dt}\Psi(t)+ C\|\na^2u\|_{L^2}^2+C+\varepsilon\|H_t\|_{L^2}^2
\ea\ee
where
\be\la{lv4.8'}\ba \Psi(t)&=\int P\div udx+\int\frac{1}{2} |H|^2 \div u dx-\int H\cdot \na u \cdot Hdx\\
&\leq C\|H\|_{L^4}^4+C\|P\|_{L^2}^2+C\|\na u\|_{L^2}^2\leq C.\ea\ee

  Moreover, it follows from \eqref{a1}$_3$  that
\be\la{lv4.9}\ba
\nu^{-1}\frac{d}{dt}\|\na H\|_{L^2}^2 +\nu^{-2}\|H_t\|_{L^2}^2+ \|\Delta H\|_{L^2}^2 &  \le C\||H| |\na u|\|_{L^2}^2+C\||u| |\na H|\|_{L^2}^2 \\
 &\le C\|H\|_{L^4}^2\|\na u\|_{L^4}^2+C\||u| |\na H|\|_{L^2}^2\\
  &\le C+\|\na^2u\|_{L^2}^2+C\||u| |\na H|\|_{L^2}^2
\ea\ee
which together with \eqref{3r1} and choosing $\varepsilon$ suitably small yields that
\be\la{lv4.10}\ba
& \frac{d}{dt}\left(\mu\|\na u\|_{L^2}^2+\nu^{-1}\|\na H\|_{L^2}^2 \right) +\|\n^{1/2}u_t\|_{L^2}^2+\nu^{-2}\|H_t\|_{L^2}^2+ \|\Delta H\|_{L^2}^2\\
  &\quad\le \frac{d}{dt}\Psi(t)+C+\|\na^2u\|_{L^2}^2+C\||u| |\na H|\|_{L^2}^2,
\ea\ee
where
\be \la{lv4.7}  \ba
\||u| |\na H|\|_{L^2}^2&=\int |u|^2\bar{x}^{-1/2}|\na H|\bar{x}^{1/2}|\na H|dx\\
&\leq C\|u \bar{x}^{-1/4}\|_{L^8}^4 \|\na H\|_{L^4}^2+C\|\na H \bar{x}^{1/2}\|_{L^2}^2\\
&\leq \frac{1}{2}\|\na^2  H\|_{L^2}^2+C\|\na H \bar{x}^{a/2}\|_{L^2}^2.
\ea\ee

Now,  integrating  \eqref{lv4.10} over $(0,T)$, along with \eqref{lv4.8'}, \eqref{lv4.7}, \eqref{pa1}, \eqref{gj10} and  \eqref{q}, we have
\be\la{gj12}
\sup_{0\le t\le T}\left(\mu\|\na u\|_{L^2}^2+\nu^{-1}\|\na H\|_{L^2}^2\right)+\int_0^T\left(\|\rho^{1/2}u_t\|_{L^2}^2+\nu^{-2}\|H_t\|_{L^2}^2+ \|\Delta H\|_{L^2}^2\right)dt\leq C.
\ee

Now, differentiating $\eqref{a1}_2$ with respect to $t$ gives
 \be\la{zb1}\ba &\n u_{tt}+\n u\cdot \na u_t-\mu\Delta u_t-( \mu+\lm)\na  \div u_t  \\ &=-\n_t(u_t+u\cdot\na u)-\n u_t\cdot\na u -\na P_t+\left(H\cdot\na H-\frac{1}{2}\na |H|^2\right)_t.\ea\ee
Multiplying \eqref{zb1} by $u_t$, then integrating  over $\O,$ we obtain after using  $\eqref{a1}_1$ that\be\ba  \la{na8}&\frac{1}{2}\frac{d}{dt} \int \n |u_t|^2dx+\int \left(\mu|\na u_t|^2+( \mu+\lm)(\div u_t)^2  \right)dx\\
 &=-2\int \n u \cdot \na  u_t\cdot u_tdx  -\int \n u \cdot\na (u\cdot\na u\cdot u_t)dx\\
  &\quad-\int \n u_t \cdot\na u \cdot  u_tdx
 +\int P_{t}\div u_{t} dx+\int\left(H\cdot\na H-\frac{1}{2}\na |H|^2\right)_t u_tdx\\
 &\triangleq \sum_{i=1}^{5}\bar{J}_i.
  \ea\ee 
Similar to the proof of  \cite[Lemma 4.3]{lx1}, we have
 \be\la{na2}\ba  \sum_{i=1}^{4}\bar{J}_i \le  \de\| \na u_{t}\|_{L^{2}}^{2}+C(\de)  \left(\| \na^{2} u \|_{L^{2}}^{2} +  \|\n^{1/2} u_{t}\|_{L^{2}}^{2}+1\right).
 \ea\ee
For the term $\bar{J}_5$,   we obtain after integration by parts  that
\be\la{na3}\ba
\bar{J}_5=-\int H_t\cdot \na u_t\cdot Hdx-\int H\cdot \na u_t\cdot H_tdx+\int H\cdot H_t \div u_tdx\triangleq \sum_{i=1}^{3} S_i.
\ea\ee


Next, differentiating $\eqref{a1}_3$ with respect to $t$ shows
  \be\la{lv4.12}\ba
H_{tt}-H_t\cdot\na u-H\cdot\na u_t+u_t\cdot\na H+u\cdot\na H_t+H_t\div u+H\div u_t=\nu \Delta H_t.
\ea\ee
 Multiplying \eqref{lv4.12} by $H_t$ and integrating the resulting equation over $\O,$ yileds that
\be\la{lv4.13}\ba
&\frac{1}{2}\frac{d}{dt} \int |H_t|^2dx+\nu \int  |\na H_t|^2dx\\
&=\int H\cdot\na u_t\cdot H_tdx-\int u_t\cdot\na H\cdot H_tdx-\int H \cdot H_t \div u_tdx\\
&\quad+\int  H_t\cdot\na u\cdot H_tdx-\int u\cdot\na H_t\cdot H_tdx-\int |H_t|^2\div u dx\triangleq \sum_{i=4}^{9} S_i.
\ea\ee

For the terms $S_i(i=1,\cdots,9)$ on the right hand of \eqref{na3} and \eqref{lv4.13}, we have
\be\la{lvbo4.14}\ba
\sum_{i=1}^{6} S_i&=-\int H_t\cdot \na u_t\cdot Hdx-\int u_t\cdot \na H\cdot H_tdx\\
&\leq  C\|H_t\|_{L^4}\|H\|_{L^4}\|\na u_t\|_{L^2}
+C\|u_t\bar{x}^{-a}\|_{L^4}\||\na H|^{1/2}\bar{x}^{a}\|_{L^4}\||\na H|^{1/2}\|_{L^4}\|H_t\|_{L^4}\\
&\leq  C\|H_t\|_{L^4}\|\na u_t\|_{L^2}
+C\left(\|\n^{1/2}u_t\|_{L^2}+\|\na u_t\|_{L^2}\right)\|\na H\bar{x}^{a/2}\|_{L^2}^{1/2}\|H_t\|_{L^4}\\
&\le \delta  \|\na u_t\|_{L^2}^2+\delta \|\n^{1/2}u_t\|_{L^2}^2+C(\delta)  \|H_t\|_{L^4}^2
+C(\delta)\|H_t\|_{L^4}^2\|\na H\bar{x}^{a/2}\|_{L^2}\\
&\le \delta  \|\na u_t\|_{L^2}^2+\delta  \|\na H_t\|_{L^2}^2+C(\delta) \|\n^{1/2}u_t\|_{L^2}^2\\
&\quad+C(\delta)  \|H_t\|_{L^2}^2+C(\delta)\|H_t\|_{L^2}^2\|\na H\bar{x}^{a/2}\|_{L^2}^2
\ea\ee
owing to \eqref{pa1} and \eqref{ljo1},
\be\la{lv4.14}\ba
\sum_{i=7}^{9} S_i& \leq C \int  |H_t|^2|\na u|dx \leq  C \|H_t\|_{L^2} \|\na H_t\|_{L^2} \|\na u\|_{L^2}\\
&\leq \delta \|\na H_t\|_{L^2}^2+C(\delta)\|H_t\|_{L^2}^2,
\ea\ee
due to  \eqref{pa1}.

 Now, adding \eqref{na8} multiplied by $t$ and \eqref{lv4.13}  multiplied by $t$ together, choosing $\delta$ suitably small and using \eqref{na2}, \eqref{na3},  \eqref{lvbo4.14}, \eqref{lv4.14}, \eqref{gj10'},  \eqref{pa1}, we have that
\be\ba\la{lv4.18}& \frac{1}{2}\frac{d}{dt} \left( t\| \n^{1/2}  u_t\|^2_{L^2}+ t\|H_t\|^2_{L^2} \right)+\mu t\| \na u_t\|^2_{L^2}+\nu t\| \na H_t\|^2_{L^2}\\
 &\le    C\left( t\|\n^{1/2} u_{t}\|_{L^{2}}^{2}+t\|H_t\|_{L^2}^2\right)+ C\left( \|\n^{1/2} u_{t}\|_{L^{2}}^{2}+\|H_t\|_{L^2}^2+1\right)
  \ea\ee
which together with Gronwall's inequality and  \eqref{gj12}  yields that \eqref{gj13}.

Finally, notice that
  \be\ba\la{AMSS11}\|\na u\|^2_{H^1}+\|\na H\|^2_{H^1}&\leq C\|\na u\|^2_{L^2}+C\|\na H\|^2_{L^2}+C\|\na^2 u\|^2_{L^2}+C\|\na^2 H\|^2_{L^2}\\
  &\leq C\|\na u\|^2_{L^2}+C\|\na H\|^2_{L^2}+C\|\na^2 u\|^2_{L^2}+C\|H_t\|^2_{L^2}\\
  &~~~~+C\||u| |\na H|\|_{L^2}^2+C\||H| |\na u|\|^2_{L^2}\\
  &\leq C\|\na u\|^2_{L^2}+C\|\na H\|^2_{L^2}+C\|\na^2 u\|^2_{L^2}\\
  &~~~~+C\|H_t\|^2_{L^2}+ \frac{1}{2}\|\na^2  H\|_{L^2}^2+C\|\na H \bar{x}^{a/2}\|_{L^2}^2
   \ea\ee
where in the second and last inequalities one has used respectively \eqref{a1}$_3$ and  \eqref{lv4.7}.
Multiplying \eqref{AMSS11}  by $t$, we obtain \eqref{gj13'} directly from \eqref{pa1}, \eqref{gj13} and \eqref{gj10'}.

  The proof of Lemma \ref{lem4.5v} is finished.

\section{\la{se4}Proofs of  Theorems  \ref{th1} and \ref{thv} }

With all the a priori estimates in Sections \ref{se3} and \ref{se5} at hand, we are
ready to prove the main result of this paper in this section.

{\it Proof of Theorem \ref{th1}.} By Lemma \ref{th0}, there exists a
$T_*>0$ such that the Cauchy problem (\ref{a1})--(\ref{h2}) has a unique strong solution $(\n,u,H)$ on $\O \times
(0,T_*]$. We will use the a priori estimates, Proposition \ref{pr1}
and Lemmas \ref{le4}-\ref{lem4.5v}, to prove the local strong
solution $(\n,u,H)$ shall exist for all time.

First,  it follows from (\ref{As1}), (\ref{As2}), (\ref{A3}), and
(\ref{co1})    that
$$ A_1(0)+A_2(0)=0,   \quad A_3(0)=0,  \quad \n_0\leq
\bar{\n}.$$   Therefore, there exists a
$T_1\in(0,T_*]$ such that (\ref{z1}) holds for $T=T_1$.

Next, set \bn \la{s1}T^*=\sup\{T\,|\,{\rm (\ref{z1}) \
holds}\}.\en Then $T^*\geq T_1>0$. Hence, for any $0<\tau<T\leq T^*$
with $T$ finite, one deduces from   \eqref{gj13} and \eqref{gj13'}
that
 \be \la{sp43}\na u,~\na H \in C([\tau ,T];L^2\cap L^q),\ee
where one has used the standard embedding
$$L^{\infty}(\tau,T; H^1)\cap H^1(\tau,T; H^{-1})\hookrightarrow C(\tau,T; L^q),\quad\texttt{for~any}~q\in[2,\infty].$$
Moreover, it follows from \eqref{pa1}, \eqref{q} and  \cite[Lemma 2.3]{L2} that \be \la{n20}\n\in C([0,T];L^1\cap H^1\cap W^{1,q}).\ee

Finally,   we claim that \be \la{s2}T^*=\infty.\ee Otherwise,
$T^*<\infty$. Then by Proposition \ref{pr1}, (\ref{z2}) holds for
$T=T^*$. It follows from  \eqref{a16}, \eqref{q},  \eqref{gj10},  \eqref{sp43} and
(\ref{n20})  that $(\n(x,T^*),u(x,T^*),H(x,T^*))$ satisfies
(\ref{co1})   except $( u(\cdot,T^*),H(\cdot,T^*))\in \dot H^\beta.$
  Thus, Lemma
\ref{th0} implies that there exists some $T^{**}>T^*$, such that
(\ref{z1}) holds for $T=T^{**}$, which contradicts (\ref{s1}).
Hence, (\ref{s2}) holds. Lemmas \ref{th0} and \ref{le4}-\ref{lem4.5v}
  thus show that $(\n,u,H)$ is in fact the unique
strong  solution defined on $\O \times(0,T]$ for any
$0<T<T^*=\infty$. Thus, the proof of Theorem  \ref{th1} is completed.

Next, we state the following well-known Gagliardo-Nirenberg inequality (see \cite{nir}) in $\mathbb{R}^3$:
\begin{lemma}
[Gagliardo-Nirenberg-3D]\la{gn-3d} For  $p\in [2,6],q\in(1,\infty), $ and
$ r\in  (3,\infty),$ there exists some generic
 constant
$C>0$ which may depend  on $p,q, $ and $r$ such that for   $f\in H^1(\mathbb{R}^3) $
and $g\in  L^q(\mathbb{R}^3 )\cap D^{1,r}(\mathbb{R}^3), $    we have \be
\la{gn1}\|f\|_{L^p(\mathbb{R}^3)}^p\le C \|f\|_{L^2(\mathbb{R}^3)}^{\frac{6-p}{2}}\|\na
f\|_{L^2(\mathbb{R}^3)}^{\frac{3p-6}{2}} ,\ee  \be
\la{gn2}\|g\|_{L^\infty\left( \mathbb{R}^3  \right)} \le C
\|g\|_{L^q(\mathbb{R}^3)}^\frac{q(r-3)}{3r+q(r-3)} \|\na g\|_{L^r(\mathbb{R}^3)}^\frac{3r}{3r+q(r-3)}.
\ee
\end{lemma}

To prove Theorem  \ref{thv}, we need the following  elementary estimates similar to those of Lemma \ref{le3} whose proof can be found in \cite[Lemma 2.2]{lxz}.
\begin{lemma} \la{lve3}
  Let $\Omega=\r^3$ and $(\rho,u,H)$ be a smooth solution of
   (\ref{a1}).
    Then there exists a generic positive
   constant $C$ depending only on $\mu$, $\lambda$  and $\nu$ such that for any $p\in [2,6]$
\be\la{hv19}
    \|{\nabla F}\|_{L^p} + \|{\nabla \o}\|_{L^p}
   \le C\left(\norm[L^p]{\rho\dot{u}}+\norm[L^p]{|H||\na H|}\right),\ee
\be\la{hv20}\ba
\norm[L^p]{F} + \norm[L^p]{\o}\le &C \left(\norm[L^2]{\rho\dot{u}}+\norm[L^2]{|H||\na H|}\right)^{(3p-6)/(2p) }\\
   &\cdot\left(\norm[L^2]{\nabla u}+ \norm[L^2]{P}+\|H\|_{L^4}^2\right)^{(6-p)/(2p)} ,
\ea\ee \be \la{hv18}
   \norm[L^p]{\nabla u} \le C \left(\norm[L^p]{F} + \norm[L^p]{\o}\right)+
   C \norm[L^p]{P }+C\||H|^2\|_{L^p},
  \ee where $F$ and $\o $ are defined in \eqref{hj1}.
\end{lemma}

{\it Proof of Theorem  \ref{thv}.} It suffices to prove \eqref{lvy8}. In fact, it follows from \cite[Proposition 3.1 and (3.10)]{lxz} that there exists some $\ve$ depending only on $\mu,\nu,\lambda,\gamma,\bar\n,\beta,$ and $M$  such that
\be \la{hv27}\ba
&\sup\limits_{1\le t<\infty}(\|\na u\|_{L^2}+\|\n\|_{L^\ga\cap L^\infty}+ \|\n^{1/2}\dot u\|_{L^2} +\|H\|_{H^2}+\|H\|_{L^3}+\|H_t\|_{L^2})\\
&\quad+\int_1^\infty(\|\na u\|_{L^2}^2+  \|\na \dot u\|_{L^2}^2+\|\n^{1/2}\dot u\|_{L^2}^2+\|\na H\|^2_{H^1}+\|H_t\|_{H^1}^2)dt\le C,
\ea\ee
provided $C_0\le \ve.$

First, we shows that
  \be \la{vo2}\sup\limits_{1\le t<\infty}\|\n\|_{L^{3/2}} \le C \ee
whose proof is completed in \cite{lx1}.

%
 Similar to \eqref{u'},  Sobolev inequality together with \eqref{u}, \eqref{hv27}  and  \eqref{vo2},  gives
\bnn \ba \|P\|_{L^2} \le&C \|(-\Delta)^{-1} \div(\n  \dot u )\|_{L^{2}}   +C\|\na u\|_{L^2} +C\|H\|_{L^4}^2 \\
\le& C\|  \n  \dot u \|_{L^{6/5}}   +C\|\na u\|_{L^2}+C\|H\|_{L^3}\| H\|_{L^6} \\
\le& C\|\n\|_{L^{3/2}}^{1/2}\| \n^{1/2}  \dot u  \|_{L^2}  +C\|\na u\|_{L^2} +C\|\na H\|_{L^2}\\
\le& C \| \n^{1/2}  \dot u  \|_{L^2} +C\|\na u\|_{L^2}+C\|\na H\|_{L^2} ,
\ea \enn
which combined with \eqref{hv27} leads to
\be \la{av16}\int_1^\infty \|P\|_{L^2}^2dt\le C. \ee

 For $p\ge 2,$ we have similarly to \eqref{vv1} that
 \be\la{vv2}\ba
\left(\| P\|_{L^p}^p \right)_t+ \frac{p\ga-1}{2\mu+\lambda}\|P\|_{L^{p+ 1}}^{p+1} &
=- \frac{p\ga-1}{2\mu+\lambda}\int P^p\left(F+\frac12|H|^2\right)dx ,\ea\ee
  which together with Holder's inequality yields
\be\la{av96}\ba
\left(\| P\|_{L^p}^p \right)_t+\frac{p\ga-1}{2(2\mu+\lambda)} \|P\|_{L^{p+1}}^{p+1}
&\le C(p) \|F\|_{L^{p+1}}^{p+1}  +C(p)\||H|^2\|_{L^{p+1}}^{p+1}.\\
\ea\ee

Now, it follows from \eqref{a1}$_3$ that
\be\la{3d1}\ba
\frac{\rm d}{{\rm d}t}\|\na H\|_{L^2}^2+  \nu^{-1}\|H_t\|_{L^2}^2+ \nu\|\na^2 H\|_{L^2}^2 \le C\|\na u\|_{L^2}^4\|\na H\|_{L^2}^2.
\ea\ee
Next, for $B(t)$  defined in \eqref{nv1}, we have by  \eqref{lbq-jia20}, \eqref{3d1}   and \eqref{hv18} that
\be\la{3d5}\ba&
(B(t)+\tilde{C}\|\na H\|_{L^2}^2)'+\|\n^{1/2}\dot{u}\|_{L^2}^{2}+\tilde{C}\nu^{-1}\|H_t\|_{L^2}^2+\tilde{C}\nu\|\na^2 H\|_{L^2}^2\\
&\leq  C\|P\|_{L^3}^3
 + C\|\nabla u\|_{L^3}^3+
C\|\na u\|_{L^2}^4\|\na H\|_{L^2}^2\\
&\leq C\|P\|_{L^3}^3  +C\|\na u\|_{L^2}^2\|\na H\|_{L^2}^2\\
&\quad+  (C\|F\|_{L^3}^3+ C\|\omega\|_{L^3}^3+ C\|P\|_{L^3}^3+ C\||H|^2\|_{L^3}^3)\\
&\le \bar{C}_1\|P\|_{L^3}^3
 + C\|F\|_{L^3}^3+ C\|\omega\|_{L^3}^3+C\|\na H\|_{L^2}^4 +
C\|\na u\|_{L^2}^4
\ea\ee
where $\tilde{C}$ is a large constant such that  $(B(t)+\tilde{C}\|\na H\|_{L^2}^2)$ satisfies \eqref{n2'}. Choosing $ \bar{C}_2\ge 1+ {2(2\mu+\lm)(\bar{C}_1+1)}/(2\ga-1)$  suitably large such that
\be \la{newn2}\ba \frac{\mu }{4}\|\nabla u\|_{L^2}^2 + \|\na H\|_{L^2}^2  + \|P\|_{L^2}^2  &\le B(t)+\tilde{C}\|\na H\|_{L^2}^2  +  \bar{C}_2\|P\|_{L^2}^2 \\
  & \le  C \|\nabla u\|_{L^2}^2+C\|\na H\|_{L^2}^2+ C \|P\|_{L^2}^2.  \ea\ee

Setting $p=2 $ in \eqref{av96},  adding \eqref{av96} multiplied by $\bar{C}_2$ to \eqref{3d5} yields that for $t\ge 1,$
\be\la{3d6}\ba
&\left(B(t)+\tilde{C}\|\na H\|_{L^2}^2+\bar{C}_2\|P\|_{L^2}^2\right)'  +\|\n^{1/2}\dot{u}\|_{L^2}^{2}+\tilde{C}\nu^{-1}\|H_t\|_{L^2}^2+\tilde{C}\nu\|\na^2 H\|_{L^2}^2+\|P\|_{L^3}^3\\
&\le C\|F\|_{L^3}^3+ C\|\omega\|_{L^3}^3 +C\|\na H\|_{L^2}^4 +C\|\na u\|_{L^2}^4
\ea\ee
owing to \eqref{hv27}. Notice that
\be\la{3d7'} \ba
\|H\|_{L^4}^{4}&\le C\|H\|_{L^3}^{2}\|H\|_{L^6}^{2} \le C\|\na H \|_{L^2}^{2}
\ea\ee
and
\be\la{3d7} \ba
\|H\cdot\na H\|_{L^2}^{2}\le \||H||\na H|\|_{L^2}^{2}\ \le C\| H \|_{L^3}^{2} \|\na H \|_{L^6}^{2}\le C \|\na^2 H \|_{L^2}^{2},
\ea\ee
which together with \eqref{hv20} gives
\be\la{3d6'}\ba
&\|F\|_{L^3}^3+ \|\omega\|_{L^3}^3  \\
&\le C\left(\|\n^{1/2}\dot{u}\|_{L^2}+\||H| |\na H|\|_{L^2}\right)^{3/2}\left(\|\na u\|_{L^2}+\|P\|_{L^2}+\|H\|_{L^4}^2\right)^{3/2}\\
&\le C\left(\|\n^{1/2}\dot{u}\|_{L^2}+\|\na^2 H \|_{L^2}\right)^{3/2}\left(\|\na u\|_{L^2}+\|P\|_{L^2}+\|\na H\|_{L^2}\right)^{3/2}\\
&\leq  \ve\|\n^{1/2}\dot{u}\|_{L^2}^{2}+\ve\|\na^2 H \|_{L^2}^2+C\|\na u\|_{L^2}^4 +C\|P\|_{L^2}^4+C\| \na H\|_{L^2}^{4}
\ea\ee
where in the last inequality one has used \eqref{hv27}.

Putting \eqref{3d6'} into \eqref{3d6} and choosing $\ve$ suitably small, then
multiplying the resulting inequality by $t$, along with Gronwall's inequality, \eqref{newn2}, \eqref{hv27} and \eqref{av16}, gives
\be\la{3d9}\ba
&\sup\limits_{1\le t<\infty}t\left(\|\na u\|_{L^2}^2+\|\na H\|_{L^2}^2+\|P\|_{L^2}^2\right)\\
&\quad+\int_1^\infty t\left(\|\n^{1/2}\dot{u}\|_{L^2}^{2}+\|H_t\|_{L^2}^2+\|\na^2 H\|_{L^2}^2+\|P\|_{L^3}^3\right)dt\leq C.
\ea\ee

Following the same arguments as  \eqref{lv3.40}, we deduce that
\be\la{3d3}\ba
\frac{\rm d}{{\rm d}t}\int\n|\dot u|^2dx +\int|\na \dot u|^2dx\le & \varepsilon\|\na \dot{u}\|_{L^2}^2+\varepsilon\|\na H_t\|_{L^2}^2+C\|\na u\|_{L^4}^4+C\|P\|_{L^4}^4\\
& +C\|\na H\|_{L^2}^4\|H_t\|_{L^2}^2+C\|\na u\|_{L^2}^2\|\na H\|_{L^2}^2\|\na^2 H\|_{L^2}^2.
\ea\ee
Furthermore, differentiating $\eqref{a1}_3$ with respect to $t$ shows
$$H_{tt}-\nu\Delta H_t=(H\cdot \na u-u\cdot \na H-H\div u)_t,$$
which together with the fact $u_t=\dot{u}-u\cdot\na u$ gives that
\be\la{3d10}\ba
&\frac{1}{2}\frac{\rm d}{{\rm d}t}\int|H_t|^2dx+\nu\int|\na H_t|^2dx\\
&\quad=\int(H_t\cdot \na u-u\cdot \na H_t-H_t\div u)H_tdx\\
&\quad\quad+\int(H \cdot \na \dot{u}-\dot{u}\cdot \na H -H \div \dot{u})H_tdx\\
&\quad\quad-\int[H \cdot \na (u\cdot\na u)- (u\cdot\na u)\cdot \na H -H \div (u\cdot\na u)]H_tdx\\
&\quad\triangleq K_1+K_2+K_3,
\ea\ee
where
\be\la{3d11}\ba
K_1+K_2&\leq C\int |H_t|^2|\na u|dx+C \|\na \dot{u}\|_{L^2}\|\na H\|_{L^2}\|H_t\|_{L^3}\\
&\leq C \|H_t\|_{L^2}^{1/2}\|\na H_t\|_{L^2}^{3/2}\|\na u\|_{L^2}+C\|\na \dot{u}\|_{L^2}\|\na H\|_{L^2}\|H_t\|_{L^2}^{1/2}\|\na H_t\|_{L^2}^{1/2}\\
&\le\varepsilon\|\na \dot{u}\|_{L^2}^2+\varepsilon\|\na H_t\|_{L^2}^2+C\|H_t\|_{L^2}^{2}(\|\na u\|_{L^2}^4+\|\na H\|_{L^2}^4)\\
K_3&\leq C \|H\|_{L^{12}}\|u\|_{L^6}\|\na u\|_{L^4}\|\na H_t\|_{L^2}\\
&\le C \| \na |H|^2\|_{L^{2}}^{1/2}\|\na u\|_{L^2}\|\na u\|_{L^4}\|\na H_t\|_{L^2}\\
 &\leq \varepsilon\|\na H_t\|_{L^2}^2+C\|\na^2 H\|_{L^2}^2\|\na u\|_{L^2}^4+C\|\na u\|_{L^4}^4. \ea\ee

Putting \eqref{3d11}  into \eqref{3d10}, adding the resulting inequality to \eqref{3d3}, and choosing $\varepsilon$ suitably small, gives
\be\la{3d14} \ba
& \frac{\rm d}{{\rm d}t}\left(\int\n|\dot u|^2dx+\int|H_t|^2dx\right)+\int|\na \dot u|^2dx+\int|\na H_t|^2dx\\
&\le C\|\na u\|_{L^4}^4+C\|P\|_{L^4}^4 +C\left(\|\na H\|_{L^2}^4+\|\na u\|_{L^2}^4\right)\left(\|H_t\|_{L^2}^2+\|\na^2 H\|_{L^2}^2\right)
\ea\ee

Next,  it follows from Gagliardo-Nirenberg inequality and \eqref{hv27} that
and
\be\la{3d19} \ba
\|\na^2H\|_{L^2}^{2}&\le C\|H_t\|_{L^2}^{2}+C\|u\cdot\na H\|_{L^2}^{2}+C\||H||\na u|\|_{L^2}^{2}\\
&\le C\|H_t\|_{L^2}^{2}+C\|\na u\|_{L^2}^{2}\|\na H\|_{L^2}\|\na^2 H\|_{L^2}\\
&\le C\|H_t\|_{L^2}^{2}+C\|\na u\|_{L^2}^{4}+C\|\na H\|_{L^2}^4+\frac{1}{2}\|\na^2 H\|_{L^2}^2.
\ea\ee
%
Moreover, by \eqref{hv20}, \eqref{hv18}, \eqref{3d7'}, \eqref{3d7} and  \eqref{3d19}, it holds that
\be\la{3d15} \ba
\|\na u\|_{L^4}^4&\le C\|F\|_{L^4}^4+C\|\omega\|_{L^4}^4+C\|P\|_{L^4}^4+C\||H|^2\|_{L^4}^4\\
&\le C\|P\|_{L^4}^4+C\left(\|\n \dot{u}\|_{L^2}+\||H| |\na H|\|_{L^2}\right)^3\left(\|\na u\|_{L^2}+\|P\|_{L^2}+\|H\|_{L^4}^2\right)\\
&\le C\|P\|_{L^4}^4+C\left(\|\n^{1/2} \dot{u}\|_{L^2}^2+\||H| |\na H|\|_{L^2}^2\right) \\
&\quad\cdot \left(\|\n^{1/2}\dot{u}\|_{L^2}^2+\||H| |\na H|\|_{L^2}^2+\|\na u\|_{L^2}^2+\|P\|_{L^2}^2+\|H\|_{L^4}^4\right)\\
&\le C\|P\|_{L^4}^4 +C\left(\|\n^{1/2}\dot{u}\|_{L^2}^2+\|H_t\|_{L^2}^2\right) \Phi(t)
 + C\left(\|\na H\|_{L^2}^4+\|\na u\|_{L^2}^4\right)\Phi(t)
\ea\ee
where
\be\la{3d15'} \ba\Phi(t)\triangleq\|\n^{1/2} \dot{u}\|_{L^2}^2+\|\na^2 H\|_{L^2}^2+\|\na u\|_{L^2}^2+\|P\|_{L^2}^2+\|\na H\|_{L^2}^2. \ea\ee
Submitting   \eqref{3d15} into \eqref{3d14}, we have
\be\la{3d17} \ba
& \frac{\rm d}{{\rm d}t}\left(\int\n|\dot u|^2dx+\int|H_t|^2dx\right)+\int|\na \dot u|^2dx+\int|\na H_t|^2dx\\
& \le \bar{C}_4\|P\|_{L^4}^4+C\left(\|\n^{1/2} \dot{u}\|_{L^2}^2+\|H_t\|_{L^2}^2\right) \Phi(t)\\
&\quad +C\left(\|\na H\|_{L^2}^4+\|\na u\|_{L^2}^4\right)\left(\Phi(t)+\|H_t\|_{L^2}^2\right)
\ea\ee
Setting $p=3$ in \eqref{av96}, and adding  \eqref{av96} multiplied by $2(2\mu+\lambda)(\bar{C}_4+1)/(3\gamma-1)$  to \eqref{3d17}, then multiplying the resulting inequality by $t^2$, lead to
\be\la{3d18} \ba
& \frac{\rm d}{{\rm d}t}\left(t^2\|\n^{1/2} \dot{u}\|_{L^2}^2+t^2\|H_t\|_{L^2}^{2}+\frac{2(2\mu+\lambda)(\bar{C}_4+1)}{3\gamma-1}t^2\|P\|_{L^3}^3\right)\\
&\qquad+t^2\|\na \dot u\|_{L^2}^2+t^2\|\na H_t\|_{L^2}^2+t^2\|P\|_{L^4}^4\\
&\le C\left(t^2\|\n^{1/2} \dot{u}\|_{L^2}^2+t^2\|H_t\|_{L^2}^{2}\right)\Phi(t)\\
&\qquad+Ct\left(\|\n^{1/2} \dot{u}\|_{L^2}^2+\|H_t\|_{L^2}^{2}+\|P\|_{L^3}^3\right)+C\left(\Phi(t)+\|H_t\|_{L^2}^2 \right)
\ea\ee
where in the last inequality one has used \eqref{3d9}.
This combined with Gronwall's inequality,  \eqref{3d9}, \eqref{hv27}  and \eqref{av16},  yields that
\be\la{lvy9}\ba   \sup\limits_{1\le t <\infty} t^2\int  \left(\n|\dot{u} |^2 +|H_t|^2+ P^3  \right)dx +  \int_{1}^\infty t^2 \left(\|\nabla\dot{u}\|^2_{L^2}+ \|P\|_{L^4}^4+\|\na H_t\|_{L^2}^2\right)dt \le C, \ea\ee
which together with  \eqref{3d7}, \eqref{3d19} and  \eqref{3d9} that
\be\la{3d21}\ba   \sup\limits_{1\le t <\infty}t^2\||H| |\na H|\|_{L^2}^{2} \le   C\sup\limits_{1\le t <\infty}t^2\|\na^2H\|_{L^2}^2\leq C . \ea\ee


Then,  \eqref{3d9}, \eqref{lvy9}  and   \eqref{3d21}  combined   with \eqref{hv18}   gives \eqref{lvy8} provided we show that for $m=1,2,\cdots,$
\be \la{lvy10}\sup\limits_{1\le t<\infty }t^{m}\|P\|_{L^{m+1}}^{m+1}+\int_0^\infty t^m\|P\|_{L^{m+2}}^{m+2}dt\le C(m).\ee
Finally, we need only to prove \eqref{lvy10}. Since \eqref{3d9} shows that \eqref{lvy10} holds for $m=1,$ we assume that \eqref{lvy10} holds for $m=n,$ that is, \be \la{lvy16}\sup\limits_{1\le t<\infty }t^{n}\|P\|_{L^{n+1}}^{n+1}+\int_1^\infty t^n\|P\|_{L^{n+2}}^{n+2}dt\le C(n) .\ee
Setting  $p=n+2 $ in \eqref{vv2} and multiplying \eqref{vv2}  by $t^{n+1} $ give
\be\la{zvo8}\ba
&\frac{2(2\mu+\lambda)}{(n+2)\ga-1}\left(t^{n+1}\| P\|_{L^{n+2}}^{n+2} \right)_t+ t^{n+1}\|P\|_{L^{n+3}}^{n+3} \\
&\quad\le C (n) t^{n}\| P\|_{L^{ n+2} }^{n+2} +C(n)   t^{n+1}\|P\|^{n+2}_{L^{n+2}} \left(\|F\|_{L^\infty}+\||H|^2\|_{L^\infty}\right). \ea\ee
It follows from  \eqref{gn1}-\eqref{gn2}, \eqref{hv19}, \eqref{3d9}, \eqref{lvy9} and \eqref{3d21}   that
\bnn \ba&\int_1^\infty \left(\|F\|_{L^\infty}+\||H|^2\|_{L^\infty}\right)dt\\
\le& C\int_1^\infty  \|\na F\|_{L^2}^{1/2} \|\na F\|_{L^6}^{1/2}dt+C\int_1^\infty  \||H| |\na H|\|_{L^2}^{1/2} \||H| |\na H|\|_{L^6}^{1/2}dt\\
\le& C\int_1^\infty  \left( \|\n^{1/2}\dot{u}\|_{L^2}+\||H| |\na H|\|_{L^2}\right)^{1/2} \left( \|\n^{1/2}\dot{u}\|_{L^6}+\||H| |\na H|\|_{L^6}\right)^{1/2} dt\\
\le& C\int_1^\infty  t^{-1/2} \left( \|\dot{u}\|_{L^6}+\|H\|_{L^\infty}\|\na H\|_{L^6}\right)^{1/2} dt\\
\le& C\int_1^\infty  t^{-1/2} \left( \|\na \dot{u}\|_{L^2}+\|H\|_{L^6}^{1/2}\|\na H\|_{L^6}^{1/2}\|\na^2 H\|_{L^2}\right)^{1/2} dt\\
\le& C\int_1^\infty  t^{-4/3}dt+ C\int_1^\infty  t^2\|\na \dot{u}\|_{L^2}^2dt+C\int_1^\infty  t^{-1/2}t^{-1/8}t^{-3/4}dt\\
\le &C
\ea\enn
which, along   with \eqref{zvo8}, \eqref{lvy16}, and Gronwall's inequality, thus   shows that \eqref{lvy10} holds for $m=n+1.$   By induction, we obtain \eqref{lvy10} and finish the proof of  \eqref{lvy8}.  The proof of Theorem  \ref{thv} is completed.

\end{document}